# MANY DIFFERENT UNIFORMITY NUMBERS OF YORIOKA IDEALS

LUKAS DANIEL KLAUSNER AND DIEGO ALEJANDRO MEJÍA

Abstract. Using a countable support product of creature forcing posets, we show that consistently, for uncountably many different functions the associated Yorioka ideals' uniformity numbers can be pairwise different. In addition we show that, in the same forcing extension, for two other types of simple cardinal characteristics parametrised by reals (localisation and anti-localisation cardinals), for uncountably many parameters the corresponding cardinals are pairwise different.

## 1. Introduction

This research forms part of the study of cardinal characteristics of the continuum which are parametrised by reals and of the forcing techniques required to separate many of them. The main motivation of this paper is to produce a forcing model where several (even uncountably many) uniformity numbers of *Yorioka ideals* are pairwise different; in doing so, we included additional types of parametrised cardinal characteristics which we refer to as *localisation* and *anti-localisation cardinals*. We first review the definition of a Yorioka ideal:

**Notation 1.1.** We fix the following terminology.

(1) For $\sigma \in (2^{<\omega})^\omega$, let $[\sigma]_\infty := \bigcap_{n<\omega} \bigcup_{i\geq n} [\sigma(i)]$, where $[t] := \{x \in 2^\omega \mid t \subseteq x\}$ for $t \in 2^{<\omega}$ (we expand this notation in Notation 1.3), and let $\mathrm{ht}_\sigma \in \omega^\omega$ be the function defined by $\mathrm{ht}_\sigma(i) := |\sigma(i)|$.

(2) For $f, g \in \omega^\omega$, we write $f \ll g$ to state "for every $m < \omega$, $f \circ \mathrm{pow}_m \leq^* g$", where $\mathrm{pow}_m \colon \omega \to \omega \colon i \mapsto i^m$.

**Definition 1.2** (Yorioka [Yor02])**.** We define Yorioka ideals in two steps:

(1) For $g \in \omega^\omega$, define $\mathcal{J}_g := \{X \subseteq 2^\omega \mid \exists\, \sigma \in (2^{<\omega})^\omega \colon X \subseteq [\sigma]_\infty \text{ and } \mathrm{ht}_\sigma = g\}$.

(2) For $f \in \omega^\omega$ increasing, define $\mathcal{I}_f := \bigcup_{g \gg f} \mathcal{J}_g$. Any family of this form is called a *Yorioka ideal.*

---

2010 *Mathematics Subject Classification.* Primary 03E17; Secondary 03E35, 03E40.

The first author was supported by the Austrian Science Fund (FWF) project P29575 "Forcing Methods: Creatures, Products and Iterations".

The second author was supported by the Austrian Science Fund (FWF) project I3081 "Filters, Ultrafilters and Connections with Forcing", the grant no. IN201711 of Dirección Operativa de Investigación – Institución Universitaria Pascual Bravo, and by the Grant-in-Aid for Early Career Scientists 18K13448, Japan Society for the Promotion of Science.

We are grateful to Martin Goldstern and Teruyuki Yorioka for giving helpful comments and suggestions to improve our paper.

Yorioka ideals are partial approximations of the $\sigma$-ideal $\mathcal{SN}$ of strong measure zero subsets of $2^\omega$. They were introduced by Yorioka [Yor02] to show that no inequality between $\mathrm{cof}(\mathcal{SN})$ and $\mathfrak{c} := 2^{\aleph_0}$ can be proved in ZFC. Though it is very easy to show that $\mathcal{SN} = \bigcap_{g\in\omega^\omega} \mathcal{J}_g$, the families of the form $\mathcal{J}_g$ are not ideals in general (Kamo and Osuga [KO08]); however, $\mathcal{I}_f$ is a $\sigma$-ideal whenever $f$ is increasing,[1] and $\mathcal{SN} = \bigcap\{\mathcal{I}_f \mid f \in \omega^\omega \text{ increasing}\}$ characterises $\mathcal{SN}$, as well. Also note that $\mathcal{I}_f \subseteq \mathcal{N}$, where $\mathcal{N}$ denotes the ideal of Lebesgue measure zero subsets of $2^\omega$.

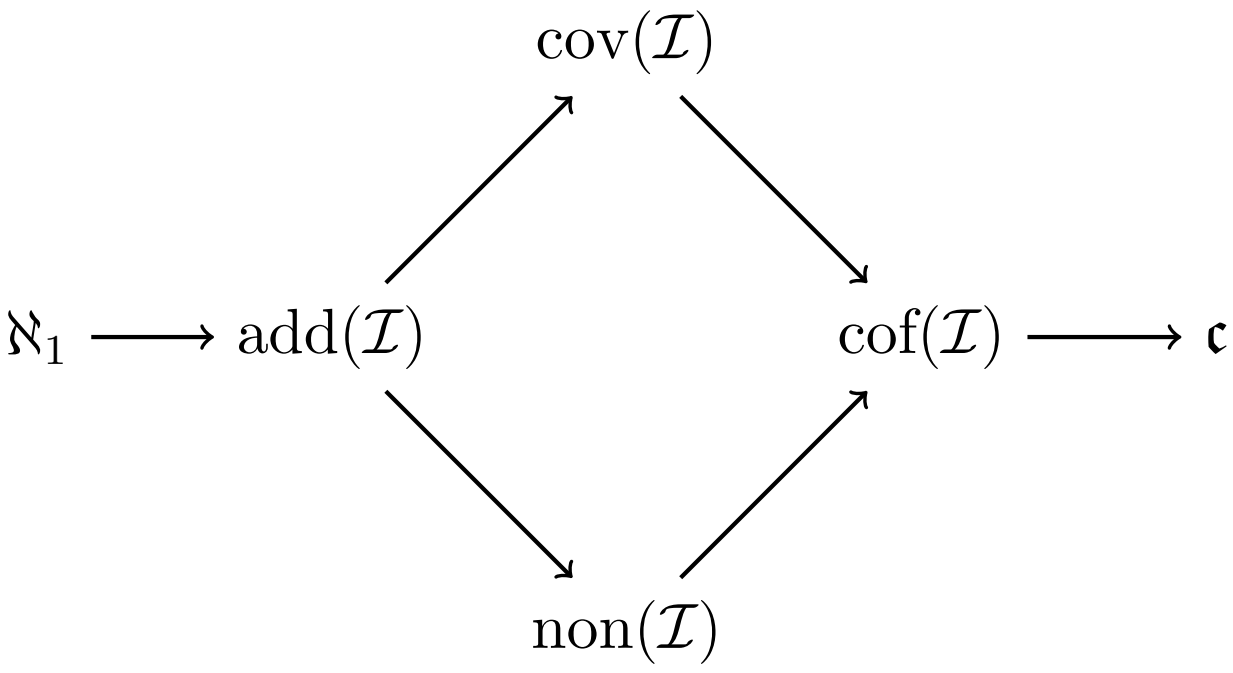

FIGURE 1. The standard inequalities for the cardinal characteristics associated with a $\sigma$-ideal $\mathcal{I}$ which contains all finite subsets and has a Borel basis.

The cardinal invariants associated with Yorioka ideals (that is, $\mathrm{add}(\mathcal{I}_f)$, $\mathrm{cov}(\mathcal{I}_f)$, $\mathrm{non}(\mathcal{I}_f)$ and $\mathrm{cof}(\mathcal{I}_f)$) have been studied for quite some time. Kamo and Osuga [KO08] showed that for any fixed increasing $f$, no other inequality consistent with the standard known inequalities (see Figure 1) can be proved in ZFC. This was improved by Cardona and the second author [CM19] by constructing a ccc poset forcing that, for some $f_0$, the four cardinals associated with $\mathcal{I}_f$ are pairwise different for all $f \geq^* f_0$. Very recently, the dependence on some $f_0$ was eliminated (see [BCM21]). On the other hand, Kamo and Osuga also showed that $\mathrm{add}(\mathcal{I}_f) \leq \mathfrak{b}$ and $\mathfrak{d} \leq \mathrm{cof}(\mathcal{I}_f)$. It is also known that $\mathrm{add}(\mathcal{N}) \leq \mathrm{add}(\mathcal{I}_f)$ and $\mathrm{cof}(\mathcal{I}_f) \leq \mathrm{cof}(\mathcal{N})$ in ZFC (attributed to Kamo; see also [CM19, section 3] for a proof), but equality cannot be proved (Cardona and the second author [CM19]).

Later, Kamo and Osuga [KO14] forced, using a ccc poset, that infinitely (even uncountably) many covering numbers of Yorioka ideals are pairwise different; moreover, continuum many pairwise different covering numbers can be forced under the assumption that a weakly inaccessible cardinal exists. The dual of this result is the main question of this paper.

**Question A.** *Is it consistent with* ZFC *that infinitely many cardinals of the form* $\mathrm{non}(\mathcal{I}_f)$ *are pairwise different?*

A key feature in Kamo and Osuga's model for infinitely many covering numbers is a relation they discovered between the covering numbers of Yorioka ideals and anti-localisation cardinals of the form $\mathfrak{c}^{\exists}_{c,h}$. The duals $\mathfrak{v}^{\exists}_{c,h}$ of these cardinals, as

[1] The proof, due to Yorioka [Yor02, Lemma 3.4], is far from trivial. An alternative proof follows from [CM19, Theorem 3.12], where it is proved that $\mathrm{add}(\mathcal{I}_f)$ is above some uncountable cardinal characteristic (without using the fact that $\mathcal{I}_f$ is an ideal itself).

well as the localisation cardinals themselves, will also play an important role in solving the question above. They are defined as follows:

**Notation 1.3.** We fix some basic notation and terminology.

(1) Given $A \subseteq \omega$ and a formula $\varphi$, we write $\forall^\infty i \in A\colon \varphi$ for $\exists n < \omega\ (\forall i \geq n, i \in A)\colon \varphi$ (i. e. all but finitely many members of $A$ satisfy $\varphi$), and $\exists^\infty i \in A\colon \varphi$ for $\forall n < \omega\ (\exists i \geq n, i \in A)\colon \varphi$ (i. e. infinitely many members of $A$ satisfy $\varphi$). We often write $\forall^\infty i$ instead of $\forall^\infty i < \omega$, and likewise for $\exists^\infty i$.
(2) Let $c = \langle c(i) \mid i < \omega\rangle$ be a sequence of non-empty sets. We write $\prod c := \prod_{i<\omega} c(i)$ and $\mathrm{seq}_{<\omega}(c) := \bigcup_{n<\omega} \prod_{i<n} c(i)$. For $t \in \mathrm{seq}_{<\omega}(c)$, let $[t] = [t]_c := \{x \in \prod c \mid t \subseteq x\}$.
(3) In addition, if $h \in \omega^\omega$, let $\mathcal{S}(c, h) := \prod_{i<\omega}[c(i)]^{\leq h(i)}$, the set of all *$h$-slaloms* (sequences of subsets of size at most $h(i)$) in $\prod c$.
(4) When $x$ and $y$ are functions with domain $\omega$, we write[2]

$$y \text{ \textit{localises} } x, \text{ denoted by } x \in^* y, \qquad \text{iff } \forall^\infty i\colon x(i) \in y(i);$$
$$y \text{ \textit{anti-localises} } x, \text{ denoted by } x \notin^\infty y, \qquad \text{iff } \forall^\infty i\colon x(i) \notin y(i).$$

(5) When, in addition, $x$ and $y$ go into the ordinal numbers, we write $x \leq y$ for $\forall i < \omega\colon x(i) \leq y(i)$ and $x \leq^* y$ for $\forall^\infty i < \omega\colon x(i) \leq y(i)$. We likewise use the notation $x < y$ and $x <^* y$.

**Definition 1.4.** Let $c = \langle c(i) \mid i < \omega\rangle$ be a sequence of non-empty sets and $h \in \omega^\omega$. We define the following cardinal characteristics:

$$\mathfrak{d}^\forall_{c,h} := \min\{|F| \mid F \subseteq \textstyle\prod c \text{ and } \neg\exists\,\varphi \in \mathcal{S}(c, h)\ \forall x \in F\colon x \in^* \varphi\}$$
$$\mathfrak{c}^\forall_{c,h} := \min\{|S| \mid S \subseteq \mathcal{S}(c, h) \text{ and } \forall x \in \textstyle\prod c\ \exists\,\varphi \in S\colon x \in^* \varphi\}$$
$$\mathfrak{d}^\exists_{c,h} := \min\{|E| \mid E \subseteq \textstyle\prod c \text{ and } \forall\,\varphi \in \mathcal{S}(c, h)\ \exists\, y \in E\colon y \notin^\infty \varphi\}$$
$$\mathfrak{c}^\exists_{c,h} := \min\{|R| \mid R \subseteq \mathcal{S}(c, h) \text{ and } \neg\exists\, y \in \textstyle\prod c\ \forall\,\varphi \in R\colon y \notin^\infty \varphi\}$$

The first two types of cardinals are referred to as *localisation cardinals*, while the latter two are referred to as *anti-localisation cardinals*.[3]

The localisation and anti-localisation cardinals are a generalisation of the cardinals used in Bartoszyński's characterisations $\mathrm{add}(\mathcal{N}) = \mathfrak{d}^\forall_{\omega,h}$ and $\mathrm{cof}(\mathcal{N}) = \mathfrak{c}^\forall_{\omega,h}$ (where $\omega$ is to be interpreted as the constant sequence $\omega$) when $h$ goes to infinity, and $\mathrm{non}(\mathcal{M}) = \mathfrak{c}^\exists_{\omega,h}$ and $\mathrm{cov}(\mathcal{M}) = \mathfrak{d}^\exists_{\omega,h}$ when $h \geq^* 1$, where $\mathcal{M}$ denotes the ideal of meagre subsets of $2^\omega$ (see e. g. [BJ95, Theorem 2.3.9, Lemmata 2.4.2 and 2.4.8] and [Bla10, Theorem 5.14 and Remark 5.15]).

Additionally, when $c$ takes infinitely many infinite values, the localisation and anti-localisation cardinals are already characterised by other well-known cardinals (see [CM19, section 3]), so these cardinals are more interesting when $c \in \omega^\omega$ and $c >^* h$. Even then, we have some trivial values: $\mathfrak{d}^\forall_{c,h}$ is finite and $\mathfrak{c}^\forall_{c,h} = \mathfrak{c}$ when

[2] The notation $x \in^\infty y$ naturally denotes $\exists^\infty i\colon x(i) \in y(i)$. We originally used $\notin^*$ instead of $\notin^\infty$, but this turned out to be unnecessarily confusing alongside our use of $\in^*$, since $\notin^\infty$ is *not* the negation of $\in^*$.

[3] A brief note on notation: $\mathfrak{c}^\forall_{c,h}$ and $\mathfrak{c}^\exists_{c,h}$ were used in [GS93, Kel08, KS09, KS12, KO14] with the meaning of *covering* $\prod c$ with slaloms, where the relations for covering are $\in^*$ and $\in^\infty$, respectively. The notation $\mathfrak{d}^\forall_{c,h}$ and $\mathfrak{d}^\exists_{c,h}$ is intended to be read as *avoiding* or *evading* such coverings.

$h$ does not go to infinity (Goldstern and Shelah [GS93]); also, $\mathfrak{c}^{\exists}_{c,h}$ is finite and $\mathfrak{v}^{\exists}_{c,h} = \mathfrak{c}$ when the quotient $\frac{h(k)}{c(k)}$ does not converge to 0 (see [CM19, section 3]). For the other cases, Figure 2 summarises the ZFC-provable inequalities relating localisation and anti-localisation cardinals to other cardinal characteristics. The additional relation between the anti-localisation cardinals and the covering and uniformity numbers of $\mathcal{N}$ was previously hinted at for the case of $h = 1$ in [GS93]; we prove it in general in Lemma 2.3.

Finally, a variation of $\mathfrak{v}^{\forall}_{\omega,h}$ appears in [BS96, Lemma 2.5] as a characterisation of the *linear evasion number* $\mathfrak{e}_\ell$ (which coincides with $\min\{\mathfrak{b}, \mathfrak{e}\}$, where $\mathfrak{e}$ denotes the *evasion number*).

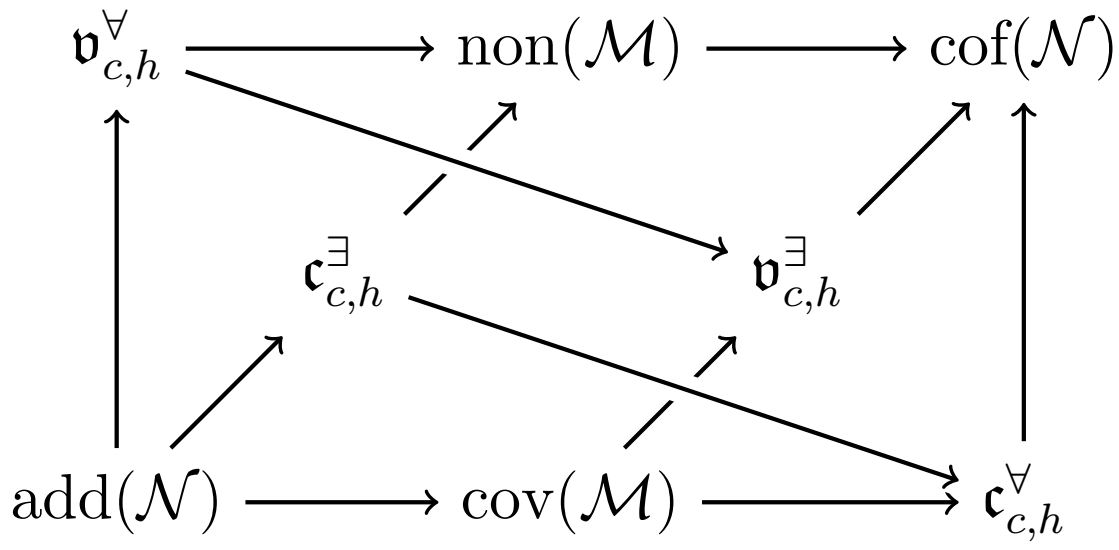

FIGURE 2. The ZFC-provable inequalities between localisation and anti-localisation cardinals and other well-known cardinal characteristics in Cichoń's diagram. Additionally, if $\sum_{i<\omega} \frac{h(i)}{c(i)} < \infty$, then $\mathrm{cov}(\mathcal{N}) \leq \mathfrak{c}^{\exists}_{c,h}$ and $\mathfrak{v}^{\exists}_{c,h} \leq \mathrm{non}(\mathcal{N})$, and conversely, if $\sum_{i<\omega} \frac{h(i)}{c(i)} = \infty$, then $\mathrm{cov}(\mathcal{N}) \leq \mathfrak{v}^{\exists}_{c,h}$ and $\mathfrak{c}^{\exists}_{c,h} \leq \mathrm{non}(\mathcal{N})$ (see Lemma 2.3).

One of the earliest appearance of these cardinals is in Miller's [Mil81] characterisations $\mathrm{non}(\mathcal{SN}) = \min\{\mathfrak{v}^{\exists}_{c,h} \mid c \in \omega^\omega\}$ ($h \geq^* 1$) and $\mathrm{add}(\mathcal{M}) = \min\{\mathfrak{b}, \mathrm{non}(\mathcal{SN})\}$.[4] Some time later, Goldstern and Shelah [GS93] proved that uncountably many cardinals of the form $\mathfrak{c}^{\forall}_{c,h}$ can be pairwise different. Kellner [Kel08] then improved this result by showing the consistency of continuum many pairwise different cardinals of the same type, and Kellner and Shelah [KS09, KS12] included similar consistency results for the type $\mathfrak{c}^{\exists}_{c,h}$; they also showed that there may be continuum many of both types of cardinals (in the same model). All of these consistency results were proved by using proper $\omega^\omega$-bounding forcing constructions with *normed creatures*.

By ccc forcing techniques, Brendle and the second author [BM14] showed that uncountably many cardinals of the form $\mathfrak{v}^{\forall}_{c,h}$ can be pairwise different – and even continuum many, under the assumption that there exists a weakly inaccessible cardinal. However, the corresponding consistency result for $\mathfrak{v}^{\exists}_{c,h}$ remained unknown.

**Question B.** *Is it consistent with* ZFC *that there are infinitely many pairwise different cardinals of the form* $\mathfrak{v}^{\exists}_{c,h}$*?*

In this paper we answer Question A and Question B in the positive. Assuming CH, we construct a forcing model where uncountably many cardinals of the form $\mathrm{non}(\mathcal{I}_f)$, $\mathfrak{v}^{\exists}_{c,h}$ and $\mathfrak{c}^{\forall}_{c,h}$ are pairwise different.

[4] By similar methods, it can also be proved that $\mathrm{cof}(\mathcal{M}) = \sup(\{\mathfrak{d}\} \cup \{\mathfrak{c}^{\exists}_{c,h} \mid c \in \omega^\omega\})$ and, whenever $h$ goes to infinity, $\mathrm{add}(\mathcal{N}) = \min(\{\mathfrak{b}\} \cup \{\mathfrak{v}^{\forall}_{c,h} \mid c \in \omega^\omega\})$ and $\mathrm{cof}(\mathcal{N}) = \sup(\{\mathfrak{d}\} \cup \{\mathfrak{c}^{\forall}_{c,h} \mid c \in \omega^\omega\})$.

**Main Theorem.** *Assume* CH. *If* $\langle \kappa_\alpha \mid \alpha \in A\rangle$ *is a sequence of infinite cardinals such that* $|A| \leq \aleph_1$ *and* $\kappa_\alpha^{\aleph_0} = \kappa_\alpha$ *for every* $\alpha \in A$*, then there is a family* $\langle (a_\alpha, d_\alpha, f_\alpha, c_\alpha, h_\alpha) \mid \alpha \in A\rangle$ *of tuples of increasing functions in* $\omega^\omega$ *and a proper* $\omega^\omega$*-bounding* $\aleph_2$*-cc poset* $\mathbb{Q}$ *that forces* $\mathfrak{v}^\exists_{c_\alpha,h_\alpha} = \mathrm{non}(\mathcal{I}_{f_\alpha}) = \mathfrak{c}^\forall_{a_\alpha,d_\alpha} = \kappa_\alpha$ *for every* $\alpha \in A$*.*

The family of increasing functions is constructed in such a way that for each $\alpha \in A$, $\mathfrak{v}^\exists_{c_\alpha,h_\alpha} \leq \mathrm{non}(\mathcal{I}_{f_\alpha}) \leq \mathfrak{c}^\forall_{a_\alpha,d_\alpha}$ is provable in ZFC. The poset $\mathbb{Q}$ is constructed by a CS (countable support) product of proper $\omega^\omega$-bounding posets, as in Goldstern and Shelah [GS93], but instead of tree-like posets we use Silver-like posets in the product. These posets are used to add generic reals that increase each $\mathfrak{v}^\exists_{c_\alpha,h_\alpha}$; in fact, forcing $\kappa_\alpha \leq \mathfrak{v}^\exists_{c_\alpha,h_\alpha}$ is not difficult. On the other hand, guaranteeing that $\mathfrak{c}^\forall_{a_\alpha,d_\alpha} \leq \kappa_\alpha$ holds in the final forcing extension requires careful definition of the increasing functions and strong combinatorics for each forcing in the product. To be more precise, each poset $\mathbb{Q}^{d_\alpha}_{c_\alpha,h_\alpha}$ (which is used to increase $\mathfrak{v}^\exists_{c_\alpha,h_\alpha}$) is obtained by a lim sup creature construction and depends on $d_\alpha$ in the sense that this function determines the *bigness* needed to not increase the cardinals $\mathfrak{c}^\forall_{a_\beta,d_\beta}$ for $\beta \neq \alpha$. Standard fusion techniques due to Baumgartner [Bau83, Bau85] are used (in section 3 and section 5) to prove that these posets and their CS products are proper with further properties (such as $\omega^\omega$-bounding with timely reading of names).

The paper is structured as follows. In section 2, we strengthen the connections Kamo and Osuga found between anti-localisation cardinals and Yorioka ideals' covering and uniformity numbers; besides, we find a (simple) connection between localisation and anti-localisation cardinals. These ensure that, in the Main Theorem, the family of increasing functions satisfies $\mathfrak{v}^\exists_{c_\alpha,h_\alpha} \leq \mathrm{non}(\mathcal{I}_\alpha) \leq \mathfrak{c}^\forall_{a_\alpha,d_\alpha}$ for every $\alpha \in A$. In section 3, for $c, h \in \omega^\omega$, we define the Silver-like poset we use to increase $\mathfrak{v}^\exists_{c,h}$ and, by incorporating a norm with sufficiently large bigness, we present sufficient conditions on functions $a, e$ which will guarantee that such a poset does not increase $\mathfrak{c}^\forall_{a,e}$. In section 4, we use the results of the previous sections to construct a family of functions as required in the Main Theorem, which is finally proved in section 5. The final section 6 is dedicated to discussions and open questions.

## 2. Yorioka Ideals and Cardinal Characteristics

For notational simplicity, we describe the cardinal characteristics we are interested in through relational systems as below.

**Definition 2.1.** A *relational system* is a triplet $\mathbf{R} := \langle X, Y, \sqsubset\rangle$ where $\sqsubset$ is a relation contained in $X \times Y$. The *cardinal characteristics associated with* $\mathbf{R}$ are

$$\mathfrak{b}(\mathbf{R}) := \min\{|B| \mid B \subseteq X \text{ and } \neg\exists\, y \in Y\ \forall x \in B\colon x \sqsubset y\},$$
$$\mathfrak{d}(\mathbf{R}) := \min\{|D| \mid D \subseteq Y \text{ and } \forall x \in X\ \exists\, y \in D\colon x \sqsubset y\}.$$

The *dual of* $\mathbf{R}$ is the relational system $\mathbf{R}^\perp := \langle Y, X, \not\sqsupset\rangle$.

Let $\mathbf{R}' := \langle X', Y', \sqsubset'\rangle$ be another relational system. A pair $(F, G)$ is a *Tukey connection from* $\mathbf{R}$ *to* $\mathbf{R}'$ if $F\colon X \to X'$, $G\colon Y' \to Y$ and for any $x \in X$ and $y' \in Y'$, $F(x) \sqsubset' y'$ implies $x \sqsubset G(y')$. When there is a Tukey connection from $\mathbf{R}$ to $\mathbf{R}'$, we say that $\mathbf{R}$ *is Tukey-below* $\mathbf{R}'$, which is denoted by $\mathbf{R} \preceq_{\mathrm{T}} \mathbf{R}'$. We say that $\mathbf{R}$ and $\mathbf{R}'$ are *Tukey-equivalent*, denoted by $\mathbf{R} \cong_{\mathrm{T}} \mathbf{R}'$, when $\mathbf{R} \preceq_{\mathrm{T}} \mathbf{R}'$ and $\mathbf{R}' \preceq_{\mathrm{T}} \mathbf{R}$.

Recall that $\mathbf{R} \preceq_{\mathrm{T}} \mathbf{R}'$ implies that $\mathfrak{d}(\mathbf{R}) \leq \mathfrak{d}(\mathbf{R}')$ and $\mathfrak{b}(\mathbf{R}') \leq \mathfrak{b}(\mathbf{R})$. Also, $\mathfrak{b}(\mathbf{R}^\perp) = \mathfrak{d}(\mathbf{R})$ and $\mathfrak{d}(\mathbf{R}^\perp) = \mathfrak{b}(\mathbf{R})$. In this section, we will use such Tukey connections to prove inequalities between cardinal invariants.

**Example 2.2.** We give two examples of relational systems.

(1) Let $\mathcal{I}$ be a family of subsets of a set $X$ that satisfies
   (i) $[X]^{<\aleph_0} \subseteq \mathcal{I}$,
   (ii) $X \notin \mathcal{I}$, and
   (iii) whenever $Y \in \mathcal{I}$ and $X \subseteq Y$, $X \in \mathcal{I}$.

   Consider the relational systems $\mathcal{I} := \langle \mathcal{I}, \mathcal{I}, \subseteq \rangle$ and $\mathbf{Cv}(\mathcal{I}) := \langle X, \mathcal{I}, \in \rangle$. Note that

$$\mathfrak{b}(\mathcal{I}) = \mathrm{add}(\mathcal{I}), \qquad \mathfrak{d}(\mathcal{I}) = \mathrm{cof}(\mathcal{I}),$$
$$\mathfrak{b}(\mathbf{Cv}(\mathcal{I})) = \mathrm{non}(\mathcal{I}), \qquad \mathfrak{d}(\mathbf{Cv}(\mathcal{I})) = \mathrm{cov}(\mathcal{I}),$$

   which are the *cardinal invariants associated with* $\mathcal{I}$. Note that $\mathbf{Cv}(\mathcal{I}) \preceq_{\mathrm{T}} \mathcal{I}$ and $\mathbf{Cv}(\mathcal{I})^\perp \preceq_{\mathrm{T}} \mathcal{I}$, so the well-known fact that $\mathrm{add}(\mathcal{I})$ is below both $\mathrm{cov}(\mathcal{I})$ and $\mathrm{non}(\mathcal{I})$ and that $\mathrm{cof}(\mathcal{I})$ is above those three is easily proved through the relational systems.

   If $\mathcal{J}$ is another family of subsets of $X$ that satisfies (i)–(iii) above and $\mathcal{I} \subseteq \mathcal{J}$, then $\mathbf{Cv}(\mathcal{J}) \preceq_{\mathrm{T}} \mathbf{Cv}(\mathcal{I})$, so $\mathrm{cov}(\mathcal{J}) \leq \mathrm{cov}(\mathcal{I})$ and $\mathrm{non}(\mathcal{I}) \leq \mathrm{non}(\mathcal{J})$.

(2) Let $c = \langle c(i) \mid i < \omega \rangle$ be a sequence of non-empty sets and $h \in \omega^\omega$. Define the relational systems $\mathbf{Lc}(c,h) := \langle \prod c, \mathcal{S}(c,h), \in^* \rangle$ and $\mathbf{aLc}(c,h) := \langle \mathcal{S}(c,h), \prod c, \not\ni^\infty \rangle$. Note that

$$\mathfrak{b}(\mathbf{Lc}(c,h)) = \mathfrak{v}^\forall_{c,h}, \qquad \mathfrak{d}(\mathbf{Lc}(c,h)) = \mathfrak{c}^\forall_{c,h},$$
$$\mathfrak{b}(\mathbf{aLc}(c,h)) = \mathfrak{c}^\exists_{c,h}, \qquad \mathfrak{d}(\mathbf{aLc}(c,h)) = \mathfrak{v}^\exists_{c,h}.$$

   If, in addition, $c' = \langle c'(i) \mid i < \omega \rangle$ is a sequence of non-empty sets, $h' \in \omega^\omega$ and $|c(i)| \leq |c'(i)|$ and $h'(i) \leq h(i)$ for all but finitely many $i < \omega$, then $\mathbf{Lc}(c,h) \preceq_{\mathrm{T}} \mathbf{Lc}(c',h')$ and $\mathbf{aLc}(c',h') \preceq_{\mathrm{T}} \mathbf{aLc}(c,h)$. Hence, $\mathfrak{v}^\forall_{c',h'} \leq \mathfrak{v}^\forall_{c,h}$, $\mathfrak{c}^\forall_{c,h} \leq \mathfrak{c}^\forall_{c',h'}$, $\mathfrak{c}^\exists_{c,h} \leq \mathfrak{c}^\exists_{c',h'}$ and $\mathfrak{v}^\exists_{c',h'} \leq \mathfrak{v}^\exists_{c,h}$.

**Lemma 2.3.** *Let $c, h \in \omega^\omega$ and assume that $c \geq 1$ and $h \geq^* 1$.*

*(a) If $\sum_{i<\omega} \frac{h(i)}{c(i)} < \infty$, then $\mathbf{Cv}(\mathcal{N}) \preceq_{\mathrm{T}} \mathbf{aLc}(c,h)^\perp$. In particular, $\mathrm{cov}(\mathcal{N}) \leq \mathfrak{c}^\exists_{c,h}$ and $\mathfrak{v}^\exists_{c,h} \leq \mathrm{non}(\mathcal{N})$.*

*(b) If $\sum_{i<\omega} \frac{h(i)}{c(i)} = \infty$, then $\mathbf{Cv}(\mathcal{N}) \preceq_{\mathrm{T}} \mathbf{aLc}(c,h)$. In particular, $\mathrm{cov}(\mathcal{N}) \leq \mathfrak{v}^\exists_{c,h}$ and $\mathfrak{c}^\exists_{c,h} \leq \mathrm{non}(\mathcal{N})$.*

Before we engage in the proof, recall that whenever $X$ is an uncountable Polish space and $\mu$ is a continuous (i. e. every singleton has measure zero) probability measure on the Borel $\sigma$-algebra of $X$, there is a Borel isomorphism $f\colon X \to 2^\omega$ that preserves the measures, i. e. $\mu(A)$ is equal to the Lebesgue measure of $f[A]$ for any Borel set $A \subseteq X$ (see e. g. [Kec95, Theorem 17.41]). Therefore, $\mathcal{N}(X) \cong_{\mathrm{T}} \mathcal{N}$ and $\mathbf{Cv}(\mathcal{N}(X)) \cong_{\mathrm{T}} \mathbf{Cv}(\mathcal{N})$ (where $\mathcal{N}(X)$ is the ideal of null sets of $X$). For the following proof, when $c \in \omega^\omega$, $c > 0$, and $\exists^\infty i\colon c(i) \geq 2$, $\prod c$ is the Polish space endowed with the product topology of the discrete spaces $c(i)$ (for $i < \omega$), and $\mu_c$ denotes the product measure of $\langle \mu_{c(i)} \mid i < \omega \rangle$, where $\mu_{c(i)}$ is the measure on the power set of $c(i)$ such that each singleton has measure $1/c(i)$. It is clear that $\mu_c$ is a continuous probability measure on the Borel $\sigma$-algebra of $\prod c$.

*Proof.* To see (a), it is enough to show that $\mathbf{Cv}(\mathcal{N}(\prod c)) \preceq_{\mathrm{T}} \mathbf{aLc}(c,h)^\perp$. Note that $F\colon \prod c \to \prod c$, defined as the identity map, and $G\colon \mathcal{S}(c,h) \to \mathcal{N}(\prod c)$, defined as $G(S) := \{x \in \prod c \mid \exists^\infty i\colon x(i) \in S(i)\}$, form the corresponding Tukey connection. We only have to check that $G(S) \in \mathcal{N}(\prod c)$. Indeed, setting $P_i^S := \{x \in \prod c \mid x(i) \in S(i)\}$, we have

$$\begin{aligned}\mu_c(G(S)) &= \lim_{n\to\infty} \mu_c\left(\bigcup_{i\geq n} P_i^S\right) \leq \lim_{n\to\infty} \sum_{i\geq n} \mu_c(P_i^S)\\ &= \lim_{n\to\infty} \sum_{i\geq n} \frac{|S(i)|}{c(i)} \leq \lim_{n\to\infty} \sum_{i\geq n} \frac{h(i)}{c(i)} = 0,\end{aligned}$$

where the last equality follows from the assumption that $\sum_{i<\omega} \frac{h(i)}{c(i)} < \infty$.

To see (b), first note that, in general, $\mathbf{aLc}(c,1) \cong_{\mathrm{T}} \mathbf{Ed}(c) := \langle \prod c, \prod c, \neq^* \rangle$, where $x \neq^* y$ means $\forall^\infty i\colon x(i) \neq y(i)$. We now first prove the following fact:

**Claim.** $\mathbf{Ed}(d) \preceq_{\mathrm{T}} \mathbf{aLc}(c,h)$, *where* $d(i) := \left\lceil \frac{c(i)}{\max\{h(i),1\}} \right\rceil$.

*Proof.* Let $h'(i) := \max\{h(i), 1\}$. Since $\forall^\infty i < \omega\colon h(i) = h'(i)$, $\mathbf{aLc}(c,h) \cong_{\mathrm{T}} \mathbf{aLc}(c,h')$ by Example 2.2 (2), so it is enough to show $\mathbf{Ed}(d) \preceq_{\mathrm{T}} \mathbf{aLc}(c,h')$. For each $i < \omega$, we can partition $c(i)$ into $d(i)$ many sets $\langle b_{i,j} \mid j < d(i) \rangle$ of size $\leq h'(i)$. Define $F\colon \prod d \to \mathcal{S}(c,h')$ by $F(x)(i) := b_{i,x(i)}$, and define $G\colon \prod c \to \prod d$ such that for any $y \in \prod c$, $G(y)(i)$ is the unique $j \in d(i)$ such that $y(i) \in b_{i,j}$. It is clear that $(F,G)$ is the required Tukey connection. □

Thanks to the preceding claim, it suffices to show $\mathbf{Cv}(\mathcal{N}(\prod d)) \preceq_{\mathrm{T}} \mathbf{Ed}(d)$. The case $d \leq^* 1$ is trivial because $\mathbf{Ed}(d) \cong_{\mathrm{T}} \mathbf{Ed}(1) \cong_{\mathrm{T}} \langle \{0\}, \{0\}, \neq \rangle$, so assume $d \not\leq^* 1$. For each $y \in \prod d$, define $G(y) := \{x \in \prod d \mid x \neq^* y\}$. Note that $\mu_d(G(y)) = \lim_{n\to\infty} \prod_{i\geq n} \left(1 - \frac{1}{d(i)}\right)$ and, using the fact that $1 + x \leq e^x$ for every real $x$,

$$\prod_{i\geq n} \left(1 - \frac{1}{d(i)}\right) \leq \prod_{i\geq n} e^{-\frac{1}{d(i)}} = e^{-\sum_{i\geq n} \frac{1}{d(i)}} = 0,$$

so $G(y) \in \mathcal{N}(\prod d)$. Defining $F\colon \prod d \to \prod d$ as the identity map, $(F,G)$ hence witnesses $\mathbf{Cv}(\mathcal{N}(\prod d)) \preceq_{\mathrm{T}} \mathbf{Ed}(d)$. □

Kamo and Osuga [KO14] proved the following connections between the covering and uniformity numbers of Yorioka ideals and anti-localisation cardinals:

(1) Let $c \in \omega^\omega$ with $c \geq^* 2$. If $g \in \omega^\omega$ and $g(n) \geq \sum_{i\leq n} \log_2 c(i)$ for all but finitely many $n < \omega$, then $\mathbf{aLc}(c,1)^\perp \preceq_{\mathrm{T}} \mathbf{Cv}(\mathcal{J}_g)$. In particular, $\mathfrak{c}^\exists_{c,1} \leq \mathrm{cov}(\mathcal{J}_g)$ and $\mathrm{non}(\mathcal{J}_g) \leq \mathfrak{v}^\exists_{c,1}$.

(2) Let $c, h \in \omega^\omega$ and $g \in \omega^\omega$ monotonically increasing. If $1 \leq^* h \leq^* c$ and $c(n) \geq 2^{g(\sum_{i\leq n} h(i) - 1)}$, then $\mathbf{Cv}(\mathcal{J}_g) \preceq_{\mathrm{T}} \mathbf{aLc}(c,h)^\perp$. In particular, $\mathrm{cov}(\mathcal{J}_g) \leq \mathfrak{c}^\exists_{c,h}$ and $\mathfrak{v}^\exists_{c,h} \leq \mathrm{non}(\mathcal{J}_g)$.

The two facts above are linked with Yorioka ideals in the sense that $\mathcal{J}_g \subseteq \mathcal{I}_f \subseteq \mathcal{J}_f$ whenever $f \ll g$ (so $\mathbf{Cv}(\mathcal{J}_f) \preceq_{\mathrm{T}} \mathbf{Cv}(\mathcal{I}_f) \preceq_{\mathrm{T}} \mathbf{Cv}(\mathcal{J}_g)$). For the purposes of this paper, we improve Kamo and Osuga's results by showing a direct connection to the Yorioka ideals without passing through a family of the form $\mathcal{J}_g$.

**Lemma 2.4.** *Let $c, h \in \omega^\omega$, let $\langle I_n \mid n < \omega\rangle$ be the interval partition such that $|I_n| = h(n)$, and let $g_{c,h} \in \omega^\omega$ be defined by $g_{c,h}(k) := \lfloor \log_2 c(n)\rfloor$ whenever $k \in I_n$. If $c \geq^* 2$, $h \geq^* 1$, $f$ is an increasing function and $g_{c,h} \gg f$, then $\mathbf{Cv}(\mathcal{I}_f) \preceq_{\mathrm{T}} \mathbf{aLc}(c,h)^\perp$. In particular, $\operatorname{cov}(\mathcal{I}_f) \leq \mathfrak{c}^\exists_{c,h}$ and $\mathfrak{v}^\exists_{c,h} \leq \operatorname{non}(\mathcal{I}_f)$.*

*Proof.* It suffices to find $F\colon 2^\omega \to \prod c$ and $G\colon \mathcal{S}(c,h) \to \mathcal{I}_f$ such that for any $y \in 2^\omega$ and $S \in \mathcal{S}(c,h)$, if $\exists^\infty n\colon F(y)(n) \in S(n)$, then $y \in G(S)$.

For each $n < \omega$, fix a one-to-one map $\iota_n\colon 2^{\lfloor \log_2 c(n)\rfloor} \to c(n)$. For $S \in \mathcal{S}(c,h)$, enumerate $S(n) =: \{m^S_{n,k} \mid k \in I_n\}$ and, whenever $k \in I_n$, define $\sigma_S(k) := \iota_n^{-1}(m^S_{n,k})$ when $m^S_{n,k} \in \operatorname{ran}\iota_n$; otherwise $\sigma_S(k)$ is allowed to be anything of length $\lfloor \log_2 c(n)\rfloor$.

Set $G(S) := [\sigma_S]_\infty$. Note that $|\sigma_S(k)| = g_{c,h}(k)$, so $\mathrm{ht}_{\sigma_s} \gg f$ and $G(S) \in \mathcal{I}_f$. On the other hand, for $y \in 2^\omega$, define $F(y)(n) := \iota_n(y{\restriction}_{\lfloor \log_2 c(n)\rfloor})$. If $F(y)(n) \in S(n)$ for infinitely many $n$, then for such $n$ there is a $k \in I_n$ such that $F(y)(n) = m^S_{n,k}$, so $\sigma_S(k) = y{\restriction}_{\lfloor \log_2 c(n)\rfloor} \subseteq y$. As there are infinitely many such $k$, we conclude that $y \in G(S)$. ☐

**Lemma 2.5.** *Let $b, g \in \omega^\omega$ such that $b \geq^* 2$, $g \geq^* 1$, let $\langle J_n \mid n < \omega\rangle$ be the interval partition such that $|J_n| = g(n)$, and let $f_{b,g} \in \omega^\omega$ be defined by $f_{b,g}(k) := \sum_{\ell\leq n}\lceil \log_2 b(\ell)\rceil$ whenever $k \in J_n$. If $f$ is an increasing function and there is some $1 \leq m < \omega$ such that $f_{b,g}(k) \leq f(k^m)$ for all but finitely many $k < \omega$, then $\mathbf{aLc}(b,g)^\perp \preceq_{\mathrm{T}} \mathbf{Cv}(\mathcal{I}_f)$. In particular, $\mathfrak{c}^\exists_{b,g} \leq \operatorname{cov}(\mathcal{I}_f)$ and $\operatorname{non}(\mathcal{I}_f) \leq \mathfrak{v}^\exists_{b,g}$.*

*Proof.* First note that for any function $g'$, $\exists\, m > 0\ \forall^\infty k\colon g'(k) \leq f(k^m)$ iff $\forall\, f' \gg f\colon g' \leq^* f'$. It suffices to define functions $F\colon \prod b \to 2^\omega$ and $G\colon \mathcal{I}_f \to \mathcal{S}(b,g)$ such that for any $y \in \prod b$ and $X \in \mathcal{I}_f$, if $F(y) \in X$, then $\exists^\infty n\colon y(n) \in G(X)(n)$.

Consider the interval partition $\langle I_n \mid n < \omega\rangle$ of $\omega$ such that $|I_n| = \lceil \log_2 b(n)\rceil$. For each $n < \omega$, fix a one-to-one function $\iota_n\colon b(n) \to 2^{I_n}$. For $y \in \prod b$, define $F(y)$ as the concatenation of the binary sequences $\langle \iota_n(y(n)) \mid n < \omega\rangle$. On the other hand, for $X \in \mathcal{I}_f$, $X \subseteq [\sigma_X]_\infty$ for some $\sigma_X \in (2^{<\omega})^\omega$ such that $\mathrm{ht}_{\sigma_X} \gg f$, so we define

$$G(X)(n) := \{\iota_n^{-1}(\sigma_X(k){\restriction}_{I_n}) \mid \sigma_X(k){\restriction}_{I_n} \in \operatorname{ran}\iota_n,\ I_n \subseteq \mathrm{ht}_{\sigma_X}(k),\ k \in J_n\}$$

(since $f_{b,h} \leq^* \mathrm{ht}_{\sigma_X}$, $\forall\, k \in J_n\colon I_n \subseteq \mathrm{ht}_{\sigma_X}(k)$ holds for all but finitely many $n$).

If $F(y) \in X$, then there is an infinite $W \subseteq \omega$ such that for any $n \in W$, $\sigma_X(k) \subseteq F(y)$ for some $k \in J_n$; so $\sigma_X(k){\restriction}_{I_n} = F(y){\restriction}_{I_n} = \iota_n(y(n))$ and hence we have $y(n) \in G(X)(n)$. ☐

The following result shows a connection between localisation cardinals and anti-localisation cardinals. This is useful to include localisation cardinals of the type $\mathfrak{c}^\forall$ in our main result.

**Lemma 2.6.** *Let $c, h \in \omega^\omega$ with $c > 0$ and $h \geq^* 1$. If $c'$ is a function with domain $\omega$ and $h' \in \omega^\omega$ such that for all but finitely many $i < \omega$, $c'(i) \geq \left|[c(i)]^{\leq h(i)} \smallsetminus \{\varnothing\}\right|$* *and* $h'(i) < {c(i)}/{h(i)}$*, then* $\mathbf{aLc}(c,h) \preceq_{\mathrm{T}} \mathbf{Lc}(c',h')$*. In particular,* $\mathfrak{v}^{\forall}_{c',h'} \leq \mathfrak{c}^{\exists}_{c,h}$ *and* $\mathfrak{v}^{\exists}_{c,h} \leq \mathfrak{c}^{\forall}_{c',h'}$*.*

*Proof.* It suffices to show the result for the case when $h(i) \geq 1$, $c'(i) = [c(i)]^{\leq h(i)} \smallsetminus \{\varnothing\}$ and $h'(i) < {c(i)}/{h(i)}$ for all $i < \omega$.[5]

Define $F\colon \mathcal{S}(c,h) \to \prod c'$ such that $F(S)(i) := S(i)$ whenever $S(i) \neq \varnothing$; otherwise $F(S)(i)$ is some arbitrary singleton. On the other hand, define $G\colon \mathcal{S}(c',h') \to \prod c$ such that $G(\varphi)(i) \in c(i) \smallsetminus \bigcup \varphi(i)$ (which is fine because this union has size $\leq h(i) \cdot h'(i) < c(i)$). It is clear that $(F,G)$ is the Tukey connection that we want. □

The next lemma is more or less the converse of the preceding one. We prove it only for completeness' sake, and it will not be used in this paper.

**Lemma 2.7.** *Let* $c, h \in \omega^\omega$ *with* $c > 0$ *and* $h \geq^* 1$*. If* $c', h' \in \omega^\omega$ *with* $c' > 0$ *such that for all but finitely many* $i < \omega$*,* $c'(i) \leq \left|[c(i)]^{\leq h(i)} \smallsetminus \{\varnothing\}\right|$ *and* $h'(i) \geq \left|[c(i)-1]^{\leq h(i)} \smallsetminus \{\varnothing\}\right|$*, then* $\mathbf{Lc}(c',h') \preceq_{\mathrm{T}} \mathbf{aLc}(c,h)$*.*

*Proof.* It suffices to show the result for the case when $h(i) \geq 1$, $c'(i) = [c(i)]^{\leq h(i)} \smallsetminus \{\varnothing\}$ and $h'(i) = \left|[c(i)-1]^{\leq h(i)} \smallsetminus \{\varnothing\}\right|$ for all $i < \omega$.

Define $F\colon \prod c' \to \mathcal{S}(c,h)$ by $F(S) := S$ and $G\colon \prod c \to \mathcal{S}(c',h')$ by $G(y)(i) := [c(i) \smallsetminus \{y(i)\}]^{\leq h(i)} \smallsetminus \{\varnothing\}$. Assume that $S \in \prod c'$, $y \in \prod c$ and $y \notin^\infty F(S)$. This means that for all but finitely many $i < \omega$, $y(i) \notin S(i)$, so $S(i) \in G(y)(i)$, that is, $S \in^* G(y)$. □

By Lemma 2.6 and Lemma 2.7 applied to $c' = c$, $h = 1$ and $h' = c - 1$, it follows that $\mathbf{aLc}(c,1) \cong_{\mathrm{T}} \mathbf{Lc}(c, c-1)$ (although it is also quite simple to prove this directly). More properties of localisation and anti-localisation cardinals can be found in e. g. [GS93, CM19].

The following lemma shows how functions should be related to get a particular chain of inequalities between anti-localisation cardinals, localisation cardinals and uniformity numbers of Yorioka ideals, which will be useful for the main result of this paper.

**Lemma 2.8.** *Assume that* $a, d, b, g, f, c, h \in \omega^\omega$ *such that*

- *(I1)* $b \geq^* 2$*,* $g \geq^* 1$ *and* $\forall^\infty i < \omega\colon {b(i)}/{g(i)} > d(i)$*,*
- *(I2)* $a(i) \geq \left|[b(i)]^{\leq g(i)} \smallsetminus \{\varnothing\}\right|$ *for all but finitely many* $i$*,*
- *(I3)* $f$ *is increasing and* $f \geq^* f_{b,g}$ *(see Lemma 2.5), and*
- *(I4)* $c \geq^* 2$*,* $h \geq^* 1$ *and* $g_{c,h} \gg f$ *(see Lemma 2.4).*

*Then* $\mathfrak{v}^{\exists}_{c,h} \leq \mathrm{non}(\mathcal{I}_f) \leq \mathfrak{v}^{\exists}_{b,g} \leq \mathfrak{c}^{\forall}_{a,d}$*.*

*Proof.* This is a direct consequence of Lemma 2.4, Lemma 2.5 and Lemma 2.6. □

[5] Note that the $c'$ we work with here and in the subsequent proof is such that $\prod c'$ already is a family of slaloms, which reduces the complexity of the Tukey connection (since we do not have to bijectively map sets to their cardinalities).

## 3. Motivational Example

We will now define a creature forcing poset (cf. [RS99] for the original seminal work on this technique) whose countable support product we will use to increase the value of $\mathfrak{v}^{\exists}_{c,h}$ to be at least some prescribed cardinal $\kappa$. Through careful choice of the function parameters (and some additional auxiliary functions), we will show in section 5 that $\mathfrak{v}^{\exists}_{c,h}$ is indeed forced to be *exactly* $\kappa$, and then use a countable support product of many such creature forcing posets with varying function parameters to force uncountably many Yorioka ideals' uniformity numbers to be different.

**Definition 3.1.** Let $c, h \in \omega^\omega$ with $c > h \geq^* 1$. We define a forcing poset $\mathbb{Q}_{c,h}$ as follows: A condition $p \in \mathbb{Q}_{c,h}$ is a sequence of *creatures* $p(n)$ such that each $p(n)$ is a non-empty subset of $[c(n)]^{\leq h(n)}$ and such that, letting the norm $\|\cdot\|_{c,h,n}$ be defined by

$$\|M\|_{c,h,n} := \max\{k \mid \forall\, Y \in [c(n)]^{\leq k}\ \exists\, X \in M \colon Y \subseteq X\},$$

$p$ fulfils $\limsup_{n\to\infty} \|p(n)\|_{c,h,n} = \infty$. The order is $q \leq p$ iff $q(n) \subseteq p(n)$ for all $n < \omega$ (i. e. stronger conditions consist of smaller subsets of $[c(n)]^{\leq h(n)}$). Note that $\mathbb{Q}_{c,h} \neq \varnothing$ iff $\limsup_{n\to\infty} h(n) = \infty$. (See Figure 3 for an illustration of the initial segment of a condition in this forcing.)

Given a condition $p$ such as above, the finite initial segments in $p{\restriction}_{k+1}$ (for $k < \omega$) are sometimes referred to as *possibilities* and denoted by $\operatorname{poss}(p, {\leq}k) := \prod_{\ell\leq k}[p(\ell)]^1 = \{\langle\{z(\ell)\} \mid \ell \leq k\rangle \mid \forall\, \ell \leq k \colon z(\ell) \in p(\ell)\}$. We may also use the notation $\operatorname{poss}(p, {<}k) := \operatorname{poss}(p, {\leq}k-1)$. When $\eta \in \operatorname{poss}(p, {\leq}k)$, we write $p \wedge \eta$ to denote $\eta^\frown p{\restriction}_{[k+1,\omega)}$.[6]

We denote the indices of the non-trivial creatures by $\operatorname{split}(p) := \langle k < \omega \mid |p(k)| > 1\rangle$ and, for $n < \omega$, denote the $n$-th member of $\operatorname{split}(p)$ by $s_n(p)$. For $p, q \in \mathbb{Q}_{c,h}$, define $q \leq_n p$ as "$q \leq p$ and $q{\restriction}_{s_n(q)+1} = p{\restriction}_{s_n(q)+1}$", which means that $q$ is stronger than $p$ and they are identical up to (including) the $n$-th non-trivial creature.

**Observation 3.2.** In this section, we often work with the set $\mathbb{Q}^*_{c,h}$ of conditions $p \in \mathbb{Q}_{c,h}$ such that $\|p(s_n(p))\|_{c,h,s_n(p)} \geq n+1$ (i. e. the $n$-th non-trivial creature has norm at least $n+1$). It is clear that $\mathbb{Q}^*_{c,h}$ is dense in $\mathbb{Q}_{c,h}$.

**Lemma 3.3.** *The poset $\mathbb{Q}_{c,h}$ adds a generic slalom $\dot{S} \in \mathcal{S}(c, h)$ such that $\mathbb{Q}_{c,h}$ forces that any ground-model real in $\prod c$ is caught infinitely often by $\dot{S}$.*

*Proof.* As the set $D_n := \{p \in \mathbb{Q}_{c,h} \mid |p(n)| = 1\}$ is dense in $\mathbb{Q}_{c,h}$, we can define $\dot{S}(n)$ as the unique member of $\bigcap_{p\in\dot{G}} p(n)$ (where $\dot{G}$ is the generic filter), which clearly is a name for a member of $[c(n)]^{\leq h(n)}$.

Now assume that $x \in V \cap \prod c$ and fix a condition $p \in \mathbb{Q}_{c,h}$ and an $n_0 < \omega$. Pick some $k \geq n_0$ such that $\|p(k)\|_{c,h,k} \geq 1$ and strengthen $p$ to $q$ by setting $q(k) = \{t\}$ for some $t \in p(k)$ that contains $x(k)$ (which exists because, by the definition of the norm, $\|M\|_{c,h,n} \geq 1$ iff $\bigcup M = c(n)$). Hence $q \leq p$ and $q \Vdash x(k) \in \dot{S}(k)$.

[6] The usual creature forcing notation (as in e. g. [FGKS17]) defines the set of possibilities more abstractly as $\operatorname{poss}(p, {\leq}k) := \prod_{\ell\leq k} p(\ell)$ and defines $p \wedge \eta$ as a condition with an extended *trunk* (a concept which we did not deem necessary to introduce in our paper). Since working with possibilities $\eta$ as sequences of singletons suffices for our proofs and is conceptually easier, we instead opted for this simpler definition.

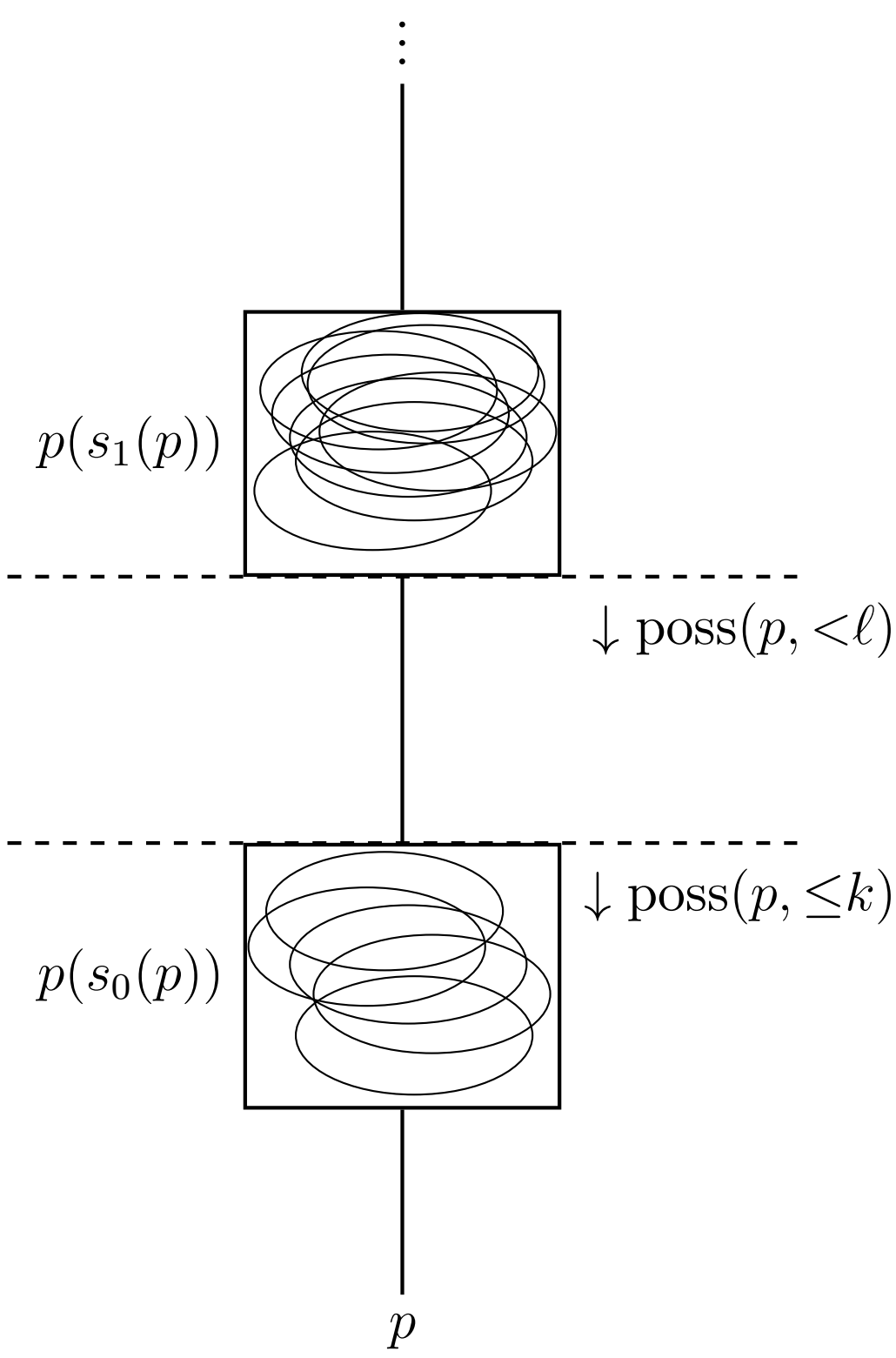

FIGURE 3. An example for the initial segment of a condition $p \in \mathbb{Q}_{c,h}$, with $k := s_0(p)$ and $\ell := s_1(p)$. Each ellipse represents one element of $[c(k)]^{\leq h(k)}$ and $[c(\ell)]^{\leq h(\ell)}$, respectively.

The previous argument shows (by density) that for any $n_0 < \omega$, $\mathbb{Q}_{c,h}$ forces $\exists\, k \geq n_0 \colon x(k) \in \dot{S}(k)$. ☐

We now show that the poset $\mathbb{Q}_{c,h}$ is indeed a proper $\omega^\omega$-bounding poset and that, moreover, $\mathbb{Q}^*_{c,h}$ satisfies strong axiom A. It is clear that $|\mathbb{Q}_{c,h}| = \mathfrak{c}$ so, assuming CH, this poset has $\aleph_2$-cc and hence preserves cardinalities and cofinalities.

**Lemma 3.4.** *If $n < \omega$, $p \in \mathbb{Q}_{c,h}$ and $D \subseteq \mathbb{Q}_{c,h}$ is open dense, then there is a condition $q \leq_n p$ in $\mathbb{Q}_{c,h}$ such that for any $\eta \in \operatorname{poss}(q, \leq s_n(q))$, $q \wedge \eta \in D$.*

*Proof.* Enumerate $\operatorname{poss}(p, \leq s_n(p)) =: \{\eta_k \mid k < m\}$. By recursion on $k \leq m$, define $p_0 := p$ and choose $p_{k+1} \leq \eta_k{}^\frown p_k{\restriction}_{[s_n(p)+1,\omega)}$ in $D$. Then

$$q := p{\restriction}_{s_n(p)+1}{}^\frown p_m{\restriction}_{[s_n(p)+1,\omega)},$$

is precisely the condition we are looking for. ☐

**Lemma 3.5.** *The poset $\mathbb{Q}^*_{c,h}$ satisfies strong axiom A. Concretely, it satisfies:*

(a) *For any $p, q \in \mathbb{Q}^*_{c,h}$ and $m \leq n < \omega$, if $p \leq_n q$, then $p \leq_m q$ and $p \leq q$.*
(b) *Whenever $\langle p_n \mid n < \omega \rangle$ is a sequence in $\mathbb{Q}^*_{c,h}$ which satisfies $p_{n+1} \leq_n p_n$ for every $n < \omega$, there is some $q \in \mathbb{Q}^*_{c,h}$ such that $q \leq_n p_n$ for every $n < \omega$.*

(c) *If* $A \subseteq \mathbb{Q}^*_{c,h}$ *is an antichain,* $p \in \mathbb{Q}^*_{c,h}$ *and* $n < \omega$*, then there is a condition* $q \leq_n p$ *in* $\mathbb{Q}^*_{c,h}$ *such that only finitely many* $r \in A$ *are compatible with* $q$*.*

*In particular,* $\mathbb{Q}_{c,h}$ *is proper and* $\omega^\omega$*-bounding.*

*Proof.* Property (a) follows immediately from the definition. To see property (b), first define $f\colon [-1,\omega) \to [-1,\omega)$ such that $f(-1) := -1$ and, for $n \geq 0$, $f(n) = s_n(p_n)$. For each $n < \omega$ define $q(k) := p_n(k)$ for any $k \in (f(n-1), f(n)]$. It is clear that $s_n(q) = f(n)$, so $q \in \mathbb{Q}^*_{c,h}$. By the definition of $q$, $q \leq_n p_n$.

We now show property (c). Set $D$ as the set of conditions $p \in \mathbb{Q}_{c,h}$ which are either stronger than some member of $A$ or incompatible with every member of $A$. It is clear that $D$ is an open dense set, so we can find $q \leq_n p$ in $\mathbb{Q}^*_{c,h}$ as in Lemma 3.4, that is, such that $q \wedge \eta \in D$ for any $\eta \in \mathrm{poss}(q, \leq s_n(q))$, which means that there is at most one $r_\eta \in A$ weaker than $q \wedge \eta$. Hence, if $r \in A$ is compatible with $q$, then it must be compatible with $q \wedge \eta$ for some $\eta \in \mathrm{poss}(q, \leq s_n(q))$ (of which there are only finitely many), so $r = r_\eta$ because $A$ is an antichain. □

A sequence $\langle p_n \mid n < \omega \rangle$ as in property (b) is usually known as a *fusion sequence*, and the $q$ obtained in its proof is known as the *fusion of* $\langle p_n \mid n < \omega \rangle$. As an application of this notion, we prove the properties of *continuous and timely reading of names* and *few possibilities* for the forcing $\mathbb{Q}_{c,h}$.

**Definition 3.6.** Let $p \in \mathbb{Q}_{c,h}$ and let $\dot\tau$ be a $\mathbb{Q}_{c,h}$-name for a function from $\omega$ into the ground model $V$.

1. We say that $p$ *reads* $\dot\tau$ *continuously* iff for each $n < \omega$ there is some $k_n$ such that for each $\eta \in \mathrm{poss}(p, \leq k_n)$, $p \wedge \eta$ decides $\dot\tau{\restriction}_n$.
2. We say that $p$ *reads* $\dot\tau$ *timely* iff for each $n \in \mathrm{split}(p)$ and each $\eta \in \mathrm{poss}(p, \leq n)$, $p \wedge \eta$ decides $\dot\tau{\restriction}_n$. (This is equivalent to continuous reading with $k_n = n$ for each $n \in \mathrm{split}(p)$.)

The term "continuous reading" refers to the fact that, in a sense, there is a continuous function from "branches" in $p$ to $V^\omega$, coded by a function from initial segments of $p$ to initial segments of $\dot\tau$. It is clear that timely reading of a name implies continuous reading thereof; we will only be using the stronger property of timely reading in the rest of this paper.

**Lemma 3.7** (timely reading of names)**.** *Let* $p \in \mathbb{Q}_{c,h}$ *and let* $\dot\tau$ *be a* $\mathbb{Q}_{c,h}$*-name for a function from* $\omega$ *into the ground model* $V$*. Then there is a condition* $q \leq p$ *in* $\mathbb{Q}_{c,h}$ *such that* $q$ *reads* $\dot\tau$ *timely.*

*Proof.* Denote by $D_n$ the set of conditions in $\mathbb{Q}_{c,h}$ which decide $\tau{\restriction}_n$, which is an open dense subset of $\mathbb{Q}_{c,h}$. By recursion on $\ell < \omega$, define an increasing function $f \in \omega^\omega$ and a sequence $\langle p_\ell \mid \ell < \omega \rangle$ of members of $\mathbb{Q}^*_{c,h}$ such that

(i) $p_0 \leq p$,
(ii) $f(\ell) = s_\ell(p_\ell)$,
(iii) $p_{\ell+1} \leq_\ell p_\ell$, and
(iv) for any $\eta \in \mathrm{poss}(p_{\ell+1}, \leq f(\ell))$, $p_{\ell+1} \wedge \eta \in D_{f(\ell)}$.

For the construction, choose any $p_0 \leq p$ in $\mathbb{Q}^*_{c,h}$ and let $f(0) := s_0(p_0)$, which clearly satisfy (i) and (ii). Now assume that $p_\ell$ and $f(\ell)$ have been constructed and that they satisfy (ii). Find $p_{\ell+1} \in \mathbb{Q}^*_{c,h}$ by application of Lemma 3.4 to $\ell$, $p_\ell$ and $D_{f(\ell)}$.

As in Lemma 3.5 (b), define $q = \langle q(k) \mid k < \omega\rangle$ such that $q{\restriction}_{(f(\ell)-1,f(\ell)]} := p_\ell{\restriction}_{(f(\ell-1),f(\ell)]}$ for any $\ell < \omega$ (letting $f(-1) := -1$). It is clear that $s_\ell(q) = f(\ell)$, even more, $\|q(f(\ell))\|_{c,h,f(\ell)} = \|p_\ell(f(\ell))\|_{c,h,f(\ell)} \geq \ell + 1$, so $q \in \mathbb{Q}^*_{c,h}$. Besides, $q \leq_\ell p_\ell$ for any $\ell < \omega$. On the other hand, by (iv), $q \wedge \eta \in D_{f(\ell)}$ for any $\eta \in \operatorname{poss}(q, {\leq} f(\ell))$. □

**Lemma 3.8** (few possibilities)**.** *Let $g \in \omega^\omega$ be a function going to infinity. If $p \in \mathbb{Q}_{c,h}$, then there is a condition $q \leq p$ in $\mathbb{Q}^*_{c,h}$ such that $|\operatorname{poss}(q, {<}s_n(q))| < g(s_n(q))$ for all $n < \omega$.*

*Proof.* We construct $f \in \omega^\omega$ and a fusion sequence $\langle p_n \mid n < \omega\rangle$ in $\mathbb{Q}^*_{c,h}$ such that

(i) $p_0 \leq p$,
(ii) $f(n) = s_n(p_n)$,
(iii) $|\operatorname{poss}(p_n, {<}f(n))| < g(f(n))$, and
(iv) $p_{n+1} \leq_n p_n$.

Start with $p_{-1} := p$ and $f(-1) := -1$. For $n < \omega$ such that $p_{n-1}$ has already been defined, choose $f(n) > f(n-1)$ in $\operatorname{split}(p_{n-1})$ such that $\|p_{n-1}(f(n))\|_{c,h,f(n)} \geq n+1$ and $|\operatorname{poss}(p_{n-1}, {\leq} f(n-1))| < g(f(n))$ (where $\operatorname{poss}(p, {\leq}{-1}) := \{\langle\ \rangle\}$). Define $p_n$ such that $p_n(k) := p_{n-1}(k)$ for every $k \in [0, f(n-1)] \cup [f(n), \omega]$ and, for $k \in (f(n-1), f(n))$, $p_n(k)$ is a singleton contained in $p_{n-1}(k)$. It is clear that $p_n \leq_{n-1} p_{n-1}$ and that $|\operatorname{poss}(p_n, {<}f(n))| = |\operatorname{poss}(p_{n-1}, {\leq}f(n-1))| < g(f(n))$.

As in previous arguments, define the fusion $q \in \mathbb{Q}^*_{c,h}$ of the sequence $\langle p_n \mid n < \omega\rangle$ by $q(k) := p_n(k)$ whenever $k \in (f(n-1), f(n)]$. It is easy to see that $q$ is as required. □

We are now interested to know for which functions $b, g \in \omega^\omega$ the poset $\mathbb{Q}_{c,h}$ will not increase the cardinal $\mathfrak{v}^\exists_{b,g}$ (i. e. any slalom in $\mathcal{S}(b,g)$ from the generic extension anti-localises some real in $\prod b$ from the ground model) or the cardinal $\mathfrak{c}^\forall_{b,g}$ (i. e. any real in $\prod b$ from the generic extension is localised by some slalom in $\mathcal{S}(b,g)$ from the ground model). This will be the key point to understand the relation between the functions so that the Main Theorem can be proved using a countable support product of our posets increasing cardinals of the form $\mathfrak{v}^\exists_{c,h}$. To ensure such *preservation*, it seems that our forcing requires the norms on the creatures to satisfy *bigness* (in the sense of [FGKS17]),[7] and for this we modify the original definition of $\mathbb{Q}_{c,h}$.

**Definition 3.9.** Let $c, h, d \in \omega^\omega$ with $c > h \geq^* 1$, $d \geq 2$, and

$$\limsup_{n\to\infty} \frac{1}{d(k)} \log_{d(k)} h(k) = \infty.$$

Recalling Definition 3.1 of $\mathbb{Q}_{c,h}$, we now define $\mathbb{Q}^d_{c,h}$ as follows: A condition $p \in \mathbb{Q}^d_{c,h}$ is a sequence of creatures $p(n)$ such that each $p(n)$ is a subset of $[c(n)]^{\leq h(n)}$ and

[7] Note that bigness is equivalent to the concept of *completeness* in the sense of [GS93].

such that, replacing $\|\cdot\|_{c,h,n}$ with the norm $\|\cdot\|^d_{c,h,n}$ defined by

$$\|M\|^d_{c,h,n} := \frac{1}{d(n)} \log_{d(n)}(\|M\|_{c,h,n} + 1),$$

$p$ fulfils $\limsup_{n\to\infty} \|p(n)\|^d_{c,h,n} = \infty$. The order is the same as on $\mathbb{Q}_{c,h}$.

Let $\mathbb{Q}^{*d}_{c,h}$ be the set of conditions in $\mathbb{Q}^d_{c,h}$ that satisfy $\|p(s_n(p))\|^d_{c,h,s_n(p)} \geq n+1$ for all $n < \omega$, which clearly is a dense subset of $\mathbb{Q}^d_{c,h}$.

The requirement $\limsup_{n\to\infty} {}^1/_{d(k)} \log_{d(k)} h(k) = \infty$ is what guarantees that the poset $\mathbb{Q}^d_{c,h}$ is non-empty. Note that $\mathbb{Q}^d_{c,h} \subseteq \mathbb{Q}_{c,h}$ and $\mathbb{Q}^{*d}_{c,h} \subseteq \mathbb{Q}^*_{c,h}$. Also note that the results proven so far remain equally valid for $\mathbb{Q}^d_{c,h}$ and $\mathbb{Q}^{*d}_{c,h}$ (the same proofs apply), in particular, that we can use the poset to increase $\mathfrak{v}^{\exists}_{c,h}$.

**Lemma 3.10** (bigness). *For any $n < \omega$, the norm $\|\cdot\|^d_{c,h,n}$ is* strongly $d(n)$-big, *i. e. whenever $M \subseteq [c(n)]^{\leq h(n)}$ and $f\colon M \to d(n)$, there is an $M' \subseteq M$ such that $f{\restriction}_{M'}$ is constant and $\|M'\|^d_{c,h,n} \geq \|M\|^d_{c,h,n} - {}^1/_{d(n)}$.*

*Proof.* This is quite similar to [FGKS17, Lemma 8.1.2 (1)]. For each $j < d(n)$, let $M_j := f^{-1}[\{j\}]$ and $m_j := \|M_j\|_{c,h,n}$. By Definition 3.1, since $m_j \leq h(n) < c(n)$, it cannot be the case that every $a \subseteq c(n)$ of size $m_j + 1$ is covered by some member of $M_j$ (otherwise $m_j + 1 \leq \|M_j\|_{c,h,n} = m_j$). Hence, there is some $a_j \subseteq c(n)$ of size $m_j + 1$ such that no member of $M_j$ contains $a_j$, thus no member of $M$ contains $a := \bigcup_{j<d(n)} a_j$, so $\|M\|_{c,h,n} \leq |a| - 1 \leq d(n) \cdot (m+1) - 1$ where $m := m_{j^*} = \max_{j<d(n)}\{m_j\}$. Therefore

$$\|M\|^d_{c,h,n} = \frac{1}{d(n)} \log_{d(n)}(\|M\|_{c,h,n} + 1) \leq \|M_{j^*}\|^d_{c,h,n} + \frac{1}{d(n)}$$

as required. □

**Corollary 3.11.** *Fix $n, m, k < \omega$. If ${}^m/_k \leq d(n)$, then, whenever $M \subseteq [c(n)]^{\leq h(n)}$ and $f\colon M \to m$, there is an $M' \subseteq M$ such that $|f[M']| \leq k$ and $\|M'\|^d_{c,h,n} \geq \|M\|^d_{c,h,n} - {}^1/_{d(n)}$.*

*Proof.* Partition $m$ into sets $\{a_j \mid j < \ell\}$ of size at most $k$ with $\ell \leq d(n)$. Define $f'\colon M \to \ell$ such that $f'(z) = j$ iff $f(z) \in a_j$. By Lemma 3.10 there is some $M' \subseteq M$ such that $f'{\restriction}_{M'}$ is constant with value $j^*$ and $\|M'\|^d_{c,h,n} \geq \|M\|^d_{c,h,n} - \frac{1}{d(n)}$. Hence $f[M'] \subseteq a_{j^*}$, so $|f[M']| \leq k$. □

Thanks to bigness, timely reading of names can be strengthened in the following way:

**Lemma 3.12** (early reading). *Let $\langle A_n \mid n < \omega\rangle$ be a sequence of non-empty finite sets, and let $\dot\tau$ be a $\mathbb{Q}^d_{c,h}$-name for a member of $\prod_{n<\omega} A_n$. If $d$ goes to infinity and $\left|\prod_{i<n} A_i\right| \leq d(n)$ for all $n < \omega$, then for any $p \in \mathbb{Q}^d_{c,h}$, there is some $q \leq p$ which* reads $\dot\tau$ early, *that is, for any $n < \omega$ and $\eta \in \operatorname{poss}(q, {<}n)$, $q \wedge \eta$ already decides $\dot\tau{\restriction}_n$.*

*Proof.* Without loss of generality, by Lemma 3.7 we may assume that $p \in \mathbb{Q}^{*d}_{c,h}$ reads $\dot\tau$ timely, that is, $p \wedge \eta$ decides $\dot\tau{\restriction}_n$ for any $\eta \in \operatorname{poss}(p, {\leq}n)$ and $n \in \operatorname{split}(p)$; by Lemma 3.8 we may also assume that $|\operatorname{poss}(p, {<}n)| < d(n)$ for any $n \in \operatorname{split}(p)$.

We construct $q(k)$ by recursion on $k < \omega$. When $k \notin \mathrm{split}(p)$, let $q(k) := p(k)$. Now assume that $k \in \mathrm{split}(p)$ and work with $p_k := q{\restriction_k}{}^\frown p{\restriction_{[k,\omega)}}$. Enumerate $\mathrm{poss}(p_k, {<}k) =: \{\eta_j \mid j < m\}$, where $m < d(k)$. By recursion on $j \leq m$, define $M_j \subseteq p(k)$ such that

(i) $M_0 = p(k)$,
(ii) $M_{j+1} \subseteq M_j$,
(iii) $\|M_{j+1}\|^d_{c,h,k} \geq \|M_j\|^d_{c,h,k} - {}^1\!/_{d(k)}$, and
(iv) there is some $r_j \in \prod_{i<k} A_i$ such that for any $t \in M_{j+1}$, $p_k \wedge (\eta_j{}^\frown\{t\})$ forces $\dot\tau{\restriction_k} = r_j$.

Assume we already have $M_j$. Define the function $h_j \colon M_j \to \prod_{i<k} A_i$ such that for any $t \in M_j$, $p_k \wedge (\eta_j{}^\frown\{t\})$ forces $\dot\tau{\restriction_k} = h_j(t)$. Hence, by Lemma 3.10, there is some $M_{j+1} \subseteq M_j$ as in (iii) such that $h_j{\restriction_{M_{j+1}}}$ is constant with value $r_j$.

Define $q(k) := M_m$. By property (iii), $\|q(k)\|^d_{c,h,k} \geq \|p(k)\|^d_{c,h,k} - {}^m\!/_{d(k)} > \|p(k)\|^d_{c,h,k} - 1$ and, by property (iv), $p_k \wedge (\eta_j{}^\frown\{t\})$ forces $\dot\tau{\restriction_k} = r_j$ for any $t \in q(k)$ and $j < m$, which means that $q{\restriction_{k+1}}{}^\frown p{\restriction_{[k+1,\omega)}} \wedge \eta_j$ forces $\dot\tau{\restriction_k} = r_j$.

By the construction, when $k \in \mathrm{split}(p)$, $\|q(k)\|^d_{c,h,k} \geq \|p(k)\|^d_{c,h,k} - 1$, so $q \in \mathbb{Q}^d_{c,h}$. If $k \in \mathrm{split}(p)$ and $\eta \in \mathrm{poss}(q, {<}k)$, then $q \wedge \eta$ decides $\dot\tau{\restriction_k}$. Now, if $k \in \omega \smallsetminus \mathrm{split}(p)$ and $\eta \in \mathrm{poss}(q, {<}k)$, then there is a unique $\eta' \in \mathrm{poss}(q, {<}\,k')$ extending $\eta$, where $k'$ is the smallest member of $\mathrm{split}(p)$ above $k$, so $q \wedge \eta = q \wedge \eta'$. As this decides $\dot\tau{\restriction_{k'}}$, it is clear that it also decides $\dot\tau{\restriction_k}$. □

The following result gives sufficient conditions on functions $a, e$ to guarantee that $\mathbb{Q}^d_{c,h}$ does not increase $\mathfrak{c}^{\forall}_{a,e}$. Hereafter, we fix the notation $c^{\nabla h}(k) := |[c(k)]^{\leq h(k)}|$.

**Lemma 3.13.** *Assume that $c, d, h \in \omega^\omega$ are as in Definition 3.9 with $d$ going to infinity, $a, e \in \omega^\omega$ with $a > 0$ and $e$ going to infinity, and that they satisfy*

*(L1) $\prod_{k<n} a(k) \leq d(n)$ and $\prod_{k<n} c^{\nabla h}(k) \leq e(n)$ for all but finitely many $n$, and*
*(L2) $\lim_{k\to\infty} \min\left\{\frac{c^{\nabla h}(k)}{e(k)}, \frac{a(k)}{d(k)}\right\} = 0$.*

*Then $\mathbb{Q}^d_{c,h}$ forces that any real in $\prod a$ is localised by some member of $\mathcal{S}(a,e) \cap V$.*

*Proof.* We show that, whenever $p \in \mathbb{Q}^d_{c,h}$ and $\dot x$ is a $\mathbb{Q}^d_{c,h}$-name for a real in $\prod a$, then there is some $\varphi \in \mathcal{S}(a,e) \cap V$ and some $q \leq p$ in $\mathbb{Q}^d_{c,h}$ forcing $\dot x(k) \in \varphi(k)$ for all but finitely many $k$. Without loss of generality, we assume that

(i) $p \in \mathbb{Q}^{*d}_{c,h}$ reads $\dot x$ early (by (L1) and Lemma 3.12),
(ii) for any $k \in \mathrm{split}(p)$, $|\,\mathrm{poss}(p, {<}k)| \leq d(k)$ and

$$\min\left\{\frac{c^{\nabla h}(k)}{e(k)}, \frac{a(k)}{d(k)}\right\} \leq \frac{1}{2\,|\,\mathrm{poss}(p, {<}k)|}$$

(by (L2) and Lemma 3.8), and
(iii) $|\,\mathrm{poss}(p, {<}k)| \leq e(k)$ for every $k < \omega$ (by (L1)).

By recursion on $k$, we construct $q(k)$ and $\varphi(k)$ according to the following case distinction: First, let $p_k := q{\restriction_k}{}^\frown p{\restriction_{[k,\omega)}}$.

When $k \notin \operatorname{split}(q)$. As $\dot{x}(k)$ is decided by the possibilities in $\operatorname{poss}(p_k, {<}k)$, the set $\varphi(k) := \{\ell \in a(k) \mid \exists\, \eta \in \operatorname{poss}(p_k, {<}k) \colon p_k \wedge \eta \Vdash \dot{x}(k) = \ell\}$ has size $\leq e(k)$ by (iii), so $\varphi(k) \in [a(k)]^{\leq e(k)}$ and $p_k$ forces $\dot{x}(k) \in \varphi(k)$. Set $q(k) := p(k)$.

When $k \in \operatorname{split}(q)$. According to (ii), we split into two subcases. If

$$|\operatorname{poss}(p_k, {<}k)| \leq \frac{e(k)}{2c^{\nabla h}(k)},$$

then $\operatorname{poss}(p_k, {\leq}k)$ has size $< e(k)$, so the set of possible values $\varphi(k)$ for $\dot{x}(k)$ has size $< e(k)$. Hence $p_k$ forces $\dot{x}(k) \in \varphi(k)$ and $\varphi(k) \in [a(k)]^{\leq e(k)}$. Set $q(k) := p(k)$.

Now consider the case when $|\operatorname{poss}(p_k, {<}k)| \leq \frac{d(k)}{2a(k)}$. Note that[8]

$$\frac{a(k)}{\left\lfloor \frac{e(k)}{|\operatorname{poss}(p_k, <k)|} \right\rfloor} \leq \frac{a(k)}{\frac{e(k)}{2 \cdot |\operatorname{poss}(p_k, <k)|}} = \frac{2a(k) \cdot |\operatorname{poss}(p_k, <k)|}{e(k)} \leq \frac{d(k)}{e(k)} \leq d(k).$$

We show how to find $q(k) \subseteq p(k)$ with $\|q(k)\|^{d(k)}_{c,h,k} \geq \|p(k)\|^{d(k)}_{c,h,k} - 1$ and $\varphi(k) \in [a(k)]^{\leq e(k)}$ such that $q{\restriction}_{(k+1)}{}^\frown p{\restriction}_{[k+1,\omega)}$ forces $\dot{x}(k) \in \varphi(k)$. Start by enumerating $\operatorname{poss}(p_k, {<}k) =: \{\eta_k \mid k < m\}$, where $m := |\operatorname{poss}(p_k, {<}k)| < d(k)$. By recursion on $j \leq m$, define $M_j \subseteq p(k)$ such that

(a) $M_0 = p(k)$,
(b) $M_{j+1} \subseteq M_j$,
(c) $\|M_{j+1}\|^d_{c,h,k} \geq \|M_j\|^d_{c,h,k} - {}^1/_{d(k)}$, and
(d) there is some $s_j \subseteq a(k)$ of size $\leq \left\lfloor \frac{e(k)}{m} \right\rfloor$ such that for any $t \in M_{j+1}$, $p_k \wedge (\eta_j{}^\frown\{t\})$ forces $\dot{x}(k) \in s_j$.

Assume we already have $M_j$. By (i), $p_k \wedge (\eta_j{}^\frown\{t\})$ decides $\dot{x}(k)$ for every $t \in M_j$, which is a value in $a(k)$. As ${}^{a(k)}/_{\lfloor \frac{e(k)}{m} \rfloor} \leq d(k)$, by Corollary 3.11 there are $M_{j+1}$ and $s_j$ as in (b)–(d). Once we have $M_m$, set $q(k) := M_m$ and $\varphi(k) := \bigcup_{j<m} s_j$.

At the end of the construction, it is clear that $q$ forces $\dot{x}(k) \in \varphi(k)$ for any sufficiently large $k$, and that $\varphi(k) \in [a(k)]^{\leq e(k)}$. ☐

**Observation 3.14.** The function $c^{\nabla h}$ is used as an upper bound of $\langle |p(k)| \mid k < \omega \rangle$ for any $p \in \mathbb{Q}^d_{c,h}$, since $p(k) \subseteq [c(k)]^{\leq h(k)}$. In the same way, we could use $c^h$ (pointwise exponentiation) instead (as long as $h \geq^* 2$), because $|[m]^{\leq k}| \leq m^k$ whenever $m, k < \omega$ and $k \neq 1$ (since $|[m]^{\leq 1}| = m + 1$). In fact, in the following sections, every instance of $c^{\nabla h}$ could be replaced by $c^h$ without affecting the results and the proofs.

Thanks to the Tukey connection constructed in Lemma 2.6, we can now easily deduce sufficient conditions on functions $b, g$ to guarantee that $\mathfrak{v}^{\exists}_{b,g}$ is not increased by $\mathbb{Q}^d_{c,h}$.

**Corollary 3.15.** *Assume that $c, d, h \in \omega^\omega$ are as in Definition 3.9 with $d$ going to infinity, $b, g, d \in \omega^\omega$ with $b, g > 0$, $e(k) := \left\lceil \frac{b(k)}{g(k)} \right\rceil - 1$ going to infinity, and that they satisfy*

*(AL1) $\prod_{k<n} b^{\nabla g}(k) \leq d(n)$ and $\prod_{k<n} c^{\nabla h}(k) \leq e(n)$ for all but finitely many $n$, and*

[8] For the first inequality, recall that $\frac{x}{2} \leq \lfloor x \rfloor$ iff $x \geq 1$.

*(AL2)* $\lim_{k\to\infty} \min\left\{\frac{c^{\nabla h}(k)}{e(k)}, \frac{b^{\nabla g}(k)}{d(k)}\right\} = 0$.

*Then $\mathbb{Q}^d_{c,h}$ forces that any slalom in $\mathcal{S}(b, g)$ is anti-localised by some member of $\prod b \cap V$.*

*Proof.* Set $a(k) := b^{\nabla g}(k)$. By Lemma 2.6, there is a definable Tukey connection $(F, G)$ which witnesses $\mathbf{aLc}(b, g) \preceq_{\mathrm{T}} \mathbf{Lc}(a, e)$ (even in forcing-generic extensions). As $a$ and $e$ satisfy the assumptions in Lemma 3.13, in any $\mathbb{Q}^d_{c,h}$-generic extension, any real in $\prod a$ is localised by some slalom in $\mathcal{S}(a, e) \cap V$. Hence any slalom in $\mathcal{S}(b, g)$ is anti-localised by some member of $\prod b \cap V$. □

## 4. Lots and Lots of Auxiliary Functions

To show that we can separate uncountably many different Yorioka ideals' uniformity numbers, we could in principle define two sequences of integers acting as universal bounds on our function parameters and auxiliary functions (similar to [GS93]). However, for the sake of clarity we will give the definition of the sequences as part of an inductive definition in the construction of the auxiliary functions. We stress that there also is an a priori definition, which would, however, make the text less readable.

Either way, we fix two sequences $n^-_k, n^+_k$ of natural numbers $\geq 2$ such that

(i) $n^-_k \cdot n^+_k < n^-_{k+1}$ for any $k < \omega$, and
(ii) $\lim_{n\to\infty} \log_{n^-_k} n^+_k = \infty$.

Given Lemma 3.13 and Corollary 3.15 above, we now make the following definition.

**Definition 4.1.** Given the bounding sequences $n^-_k, n^+_k$, we call a family $\mathcal{F} = \langle (a_\alpha, d_\alpha, b_\alpha, g_\alpha, f_\alpha, c_\alpha, h_\alpha) \mid \alpha \in A \rangle$ of tuples of increasing functions in $\omega^\omega$ *suitable* if it fulfils the following properties for any $\alpha \in A$:

(S1) The functions $a_\alpha$, $d_\alpha$, $b_\alpha$, $g_\alpha$, $b_\alpha^{\nabla g_\alpha}$, ${}^{b_\alpha}\!/\!{}_{g_\alpha}$, $h_\alpha$ and $c_\alpha^{\nabla h_\alpha}$ are bounded from below by $n^-_k$ and bounded from above by $n^+_k$, i. e. for any $k < \omega$, we have

$$a_\alpha(k), d_\alpha(k), b_\alpha(k), g_\alpha(k), b_\alpha^{\nabla g_\alpha}(k), \frac{b_\alpha(k)}{g_\alpha(k)}, h_\alpha(k), c_\alpha^{\nabla h_\alpha}(k) \in [n^-_k, n^+_k].$$

(S2) $h_\alpha < c_\alpha$ and $\limsup_{k\to\infty} \frac{1}{d_\alpha(k)} \log_{d_\alpha(k)}(h_\alpha(k) + 1) = \infty$.
(S3) ${}^{b_\alpha}\!/\!{}_{g_\alpha} > d_\alpha$.
(S4) $a_\alpha \geq b_\alpha^{\nabla g_\alpha}$.
(S5) There is some $\ell > 0$ such that $f_{b_\alpha,g_\alpha} \leq^* f_\alpha \circ \mathrm{pow}_\ell$ for $f_{b_\alpha,g_\alpha}$ as in Lemma 2.5.[9]
(S6) $f_\alpha \ll g_{c_\alpha,h_\alpha}$ for $g_{c_\alpha,h_\alpha}$ as in Lemma 2.4.
(S7) For any $\beta \in A$, if $\beta \neq \alpha$, then

$$\lim_{k\to\infty} \min\left\{\frac{c_\beta^{\nabla h_\beta}(k)}{d_\alpha(k)}, \frac{a_\alpha(k)}{d_\beta(k)}\right\} = 0.$$

[9] We can simplify property (S5) by restricting it to $\ell = 1$: Assume the statement holds for some $\ell' > 1$. Let $f'_\alpha = f_\alpha \circ \mathrm{pow}_{\ell'}$. Then $f_{b_\alpha,g_\alpha} \leq^* f'_\alpha$ and $f'_\alpha \ll g_{c_\alpha,h_\alpha}$ still holds (by the definition of $\ll$). This is why we can work with $f'_\alpha$ instead of $f_\alpha$, in which case property (S5) is already satisfied for $\ell = 1$ and property (S6) remains true for $f'_\alpha$.

Properties (S3)–(S6) ensure that, by Lemma 2.8, $\mathfrak{v}^{\exists}_{c_\alpha,h_\alpha} \leq \mathrm{non}(\mathcal{I}_{f_\alpha}) \leq \mathfrak{v}^{\exists}_{b_\alpha,g_\alpha} \leq \mathfrak{c}^{\forall}_{a_\alpha,d_\alpha}$. When we prove the Main Theorem in the following section, we will aim to force all these cardinals to be equal to some predetermined $\kappa_\alpha$ for each $\alpha$. We therefore increase the cardinals $\mathfrak{v}^{\exists}_{c_\alpha,h_\alpha}$ by using a countable support product of posets of the form $\mathbb{Q}^{d_\alpha}_{c_\alpha,h_\alpha}$, while at the same time ensuring that $\mathfrak{c}^{\forall}_{a_\alpha,d_\alpha}$ does not exceed the desired value. Property (S2) guarantees that each poset $\mathbb{Q}^{d_\alpha}_{c_\alpha,h_\alpha}$ is non-empty and, thanks to Lemma 3.13, property (S7) will ensure that $\mathbb{Q}^{d_\beta}_{c_\beta,h_\beta}$ will not increase $\mathfrak{c}^{\forall}_{a_\alpha,d_\alpha}$ for $\beta \neq \alpha$. To this end, it is also necessary that $\prod_{\ell<k} a_\alpha(\ell) \leq d_\beta(k)$ for all but finitely many $k$ (property (L1) in Lemma 3.13), which is a consequence of property (S1) because $\prod_{\ell<k} a_\alpha(\ell) \leq \prod_{\ell<k} n^+_\ell < n^-_k \leq d_\beta(k)$ (the second inequality follows from property (i) above Definition 4.1 by induction; note that the case $\alpha = \beta$ is also allowed). In fact, property (S1) is what allows us to explicitly control the functions and the number of possibilities of the creatures in the poset.

**Theorem 4.2.** *There are bounding sequences $n^-_k, n^+_k$ such that there is an uncountable suitable family $\mathcal{F} = \langle (a_\alpha, d_\alpha, b_\alpha, g_\alpha, f_\alpha, c_\alpha, h_\alpha) \mid \alpha \in 2^\omega \rangle$.*

*Proof.* For clarity, we first explain how to construct a suitable family for one single $\alpha$, that is, a suitable tuple of the form $(a, d, b, g, f, c, h)$ (as in Definition 4.1 for the case $|A| = 1$, for which (S7) is vacuous). For motivational purposes, assume we have already defined $n^-_k$ and $d(k) \geq n^-_k$, and the values of all the functions at $\ell < k$ (with the exception of $f$, which depends on some interval partition as indicated in the following argument). From here, the order in which the values should be defined is $h(k)$, $g(k)$, $b(k)$, $f{\restriction}_{I_k}$ (where $I_k$ depends on the values of $g$ up to $k$), $c(k)$, $a(k)$ and $n^+_k$. For $k = 0$, $n^-_0$ and $d(0)$ can be chosen arbitrarily such that $2 < n^-_0 < d(0)$.

To have property (S2), it suffices to define $h(k)$ such that $\frac{1}{d(k)} \log_{d(k)}(h(k)+1) \geq k+1$, so we let $h(k) := d(k)^{(k+1)\cdot d(k)}$. (Later, when we define $c(k)$, it will be clear that $c(k) > h(k)$.) The value of $g(k)$ is determined so that $f(k)$ (which will be defined a few steps later) can satisfy property (S6), but for now we will just assume that we have already defined it and will later explain how to define it according to our needs. (The value of $g(k)$ will only depend on $h(k)$ as well as $h(\ell)$ and $g(\ell)$ for $\ell < k$.)

Now, to fulfil (S3), we let $b(k) := 2^{g(k)+d(k)}$ (such that $\log_2 b(k)$ will be an integer, making the definition of $f_{b,g}$ a bit nicer). Now recall that $f_{b,g}$ is defined, along the interval partition $\langle I_n \mid n < \omega \rangle$ of $\omega$ which satisfies $|I_n| = g(n)$, as $f_{b,g}(j) = \sum_{\ell \leq n} \log_2 b(\ell)$ for any $j \in I_n$. In our construction, we may assume that we have already defined $g(\ell)$ and $I_\ell$ for any $\ell \leq k$, so we explain how to define $f{\restriction}_{I_k}$ such that property (S5) is satisfied – even more, such that $f_{b,g}(j) \leq f(j)$ for any $j \in I_k$. We also need to make sure that $f$ is increasing, so we define

$$f(j) := \sum_{\ell \leq k} \log_2 b(\ell) + j - \min(I_k)$$

for any $j \in I_k$. This definition allows $f$ to be increasing when attached to the previous intervals because

$$\begin{aligned} f(\max(I_{k-1})) &= \sum_{\ell\leq k-1} \log_2 b(\ell) + \max(I_{k-1}) - \min(I_{k-1}) \\ &= \sum_{\ell\leq k-1} \log_2 b(\ell) + g(k-1) - 1 < \sum_{\ell\leq k} \log_2 b(\ell) = f(\min(I_k)). \end{aligned}$$

Recall that $g_{c,h}$ is defined, along the interval partition $\langle J_n \mid n < \omega\rangle$ of $\omega$ which satisfies $|J_n| = h(n)$, as $g_{c,h}(j) = \lfloor \log_2 c(n)\rfloor$ for all $j \in J_n$. In order to have property (S6), it suffices to define $c(k)$ such that $f(j^{k+2}) \leq \log_2 c(k)$ for any $j \in J_k$. If we can ensure that $j^{k+2} \in I_k$ for any $j \in J_k$, then

$$\begin{aligned} f(j^{k+2}) &= \sum_{\ell\leq k} \log_2 b(\ell) + j^{k+2} - \min(I_k) \\ &\leq \sum_{\ell\leq k} \log_2 b(\ell) + \max(I_k) - \min(I_k) < \sum_{\ell\leq k} \log_2 b(\ell) + g(k), \end{aligned}$$

so it suffices to define

$$c(k) := 2^{\sum_{\ell\leq k}\log_2 b(l)+g(k)} = 2^{g(k)}\cdot \prod_{\ell\leq k} b(\ell).$$

Now, to ensure that $j^{k+2} \in I_k$ for any $j \in J_k$, it suffices to have that $\min(J_\ell)^{\ell+1} = \min(I_\ell)$ for any $\ell \leq k+1$ (as $\min(J_0) = 0$ and $\min(J_{\ell+1}) = \min(J_\ell) + h(\ell)$ as well as $\min(I_{\ell+1}) = \min(I_\ell) + g(\ell)$, the values of $\min(J_{k+1})$ and $\min(I_{k+1})$ are already known). The choice of $g(k)$, which we had postponed until now, will guarantee this fact. So we assume that we had already ensured $\min(J_k)^{k+1} = \min(I_k)$ from the beginning (so, before $n_k^-$ is even determined; also note that this is true for $k = 0$) and, after defining $h(k)$, we let $g(k) := \min(J_{k+1})^{k+2} - \min(I_k)$. This implies that $\min(I_{k+1}) = \min(I_k) + g(k) = \min(J_{k+1})^{k+2}$ as required.

We can choose $a(k)$ above $c^{\triangledown h}(k)$ and $b^{\triangledown g}(k)$, and define $n_k^+ := a(k)$, which guarantees property (S4). Note that by our definitions, we have

$$n_k^- < d(k) < h(k) < g(k) < b(k) < c(k) < c^{\triangledown h}(k) < a(k) = n_k^+$$

as well as

$$n_k^- < d(k) < \frac{b(k)}{g(k)} < b(k) < b^{\triangledown g}(k) < a(k) = n_k^+,$$

which together ensure property (S1).

We can then choose $n_{k+1}^-$ above $n_k^- \cdot n_k^+$ (note that $\log_{n_k^-} n_k^+ \geq \log_{d(k)}(h(k)+1) > k$) and define $d(k+1)$ to be larger than $n_{k+1}^-$, and then the iterative construction continues for $k+1$.

Now, we define both the bounding sequences and the suitable family $\mathcal{F}$ of size continuum inductively. More or less the same strategy as above will suffice to define the values at a fixed $k$, but some extra work is needed in order to guarantee property (S7). We construct tuples of functions $\langle (a_t, d_t, b_t, g_t, f_t, c_t, h_t) \mid t \in 2^{<\omega}\rangle$ such that

(i) $|a_t| = |d_t| = |b_t| = |g_t| = |c_t| = |h_t| = |t|$ and $|f_t| = \sum_{k<|t|} g_t(k)$, and
(ii) if $t \subseteq t'$ in $2^{<\omega}$, then $a_t \subseteq a_{t'}$, and likewise for the other functions.

Using these, we define $a_\alpha := \bigcup_{n<\omega} a_{\alpha\restriction_n}$ for $\alpha \in 2^\omega$, and likewise for the other functions, as well as $n_k^- := d_{\bar{0}}(k) - 1$ and $n_k^+ := a_{\bar{1}}(k)$ (where $\bar{e} = \langle e, e, e, \dots \rangle$ for $e \in \{0, 1\}$) and claim that these are as desired in the theorem.

Denote by $\lhd$ the lexicographic order in both $2^\omega$ and $2^n$ for any $n < \omega$. We construct $\langle (a_t, d_t, b_t, g_t, f_t, c_t, h_t) \mid t \in 2^n \rangle$ by recursion on $n$. When $n = 0$, it is clear that $a_{\langle\,\rangle}$, $d_{\langle\,\rangle}$, etc. are defined as the empty sequence. Assuming we have reached step $n$ of the construction, we show how to obtain the functions for $t \in 2^{n+1}$. We do this in the following two steps:

1. Assuming we have already defined $d_t(n)$, we define $h_t(n)$, $g_t(n)$ and so on in the same way as in the definitions for the simple case $|A| = 1$.
2. By recursion on $\langle 2^{n+1}, \lhd \rangle$, we construct $d_t(n)$ for $t \in 2^{n+1}$.

We first show step 1. At this point, we have sequences $\langle I_\ell^{t\restriction_{\ell+1}} \mid \ell < n \rangle$ and $\langle J_\ell^{t\restriction_{\ell+1}} \mid \ell < n \rangle$ of consecutive intervals covering an initial segment of $\omega$ such that $|I_\ell^{t\restriction_{\ell+1}}| = g_t(\ell)$ and $|J_\ell^{t\restriction_{\ell+1}}| = h_t(\ell)$. In fact, $I_\ell^{t\restriction_{\ell+1}} = [|f_{t\restriction_\ell}|, |f_{t\restriction_\ell}| + g_t(\ell))$ for any $\ell < n$. Note that we can already define $\min(I_n^t)$ and $\min(J_n^t)$ at this stage (considering that we are constructing interval partitions of $\omega$). Further assume that $\min(I_n^t) = \min(J_n^t)^{n+1}$, which is trivially true for $n = 0$.

Define $h_t(n) := d_t(n)^{(n+1)\cdot d_t(n)}$, $g_t(n) := (\max(J_n^t) + 1)^{n+2} - \min(I_n^t)$, where $J_n^t := [\min(J_n^t), \min(J_n^t) + h_t(n))$, and $b_t(n) := 2^{g_t(n) + d_t(n)}$. Let $I_n^t := [|f_{t\restriction_n}|, |f_{t\restriction_n}| + g_t(n))$ and, for $k \in I_n^t$, define

$$f_t(k) := \sum_{\ell \leq n} \log_2 b_t(\ell) + k - |f_{t\restriction_n}|.$$

By definition, $|f_t| = |f_{t\restriction_n}| + g_t(n) = \max(I_n^t) + 1$. Finally, let

$$c_t(n) := 2^{g_t(n)} \cdot \prod_{\ell \leq n} b_t(\ell)$$

and choose $a_t(n)$ above $c_t^{\triangledown h_t}(n)$ and $b_t^{\triangledown g_t}(n)$. As in the case of $|A| = 1$, $d_t(n) < h_t(n) < g_t(n) < b_t(n) < c_t(n) < a_t(n)$ (it will be clear from the next step that $d_t$ will always take values $\geq 2$) and the other three functions to be bounded lie somewhere in between.

To see step 2, start by choosing $d_{\bar{0}\restriction_{n+1}}(n) > d_{\bar{0}\restriction_n}(n-1) \cdot a_{\bar{1}\restriction_n}(n-1) + 2$ (in case $n = 0$, just choose $d_{\bar{0}\restriction_1}(0)$ to be any number above 2). Now assume that we have $d_t(n)$ and that $t^+ \in 2^{n+1}$ is the immediate successor of $t$ with respect to $\lhd$. Let $d_{t^+}(n) := (n+1) \cdot a_t(n)$. This completes the construction. Note that, whenever $t, t' \in 2^{n+1}$ and $t \lhd t'$, then $\frac{a_t(n)}{d_{t'}(n)} \leq \frac{1}{n+1}$, which is what ultimately guarantees property (S7).

Now we finally show that $\mathcal{F} = \langle (a_\alpha, d_\alpha, b_\alpha, g_\alpha, f_\alpha, c_\alpha, h_\alpha) \mid \alpha \in 2^\omega \rangle$ is a suitable family for the bounding sequences $n_k^-$ and $n_k^+$, that is, it satisfies the properties in Definition 4.1. Properties (S1)–(S6) are immediate by construction, as in the case $|A| = 1$. Note that $I_n^\alpha := I_n^{\alpha\restriction_{n+1}}$ and $J_n^\alpha := J_n^{\alpha\restriction_{n+1}}$ (for $n < \omega$) define interval partitions of $\omega$ such that $|I_n^\alpha| = g_\alpha(n)$ and $|J_n^\alpha| = h_\alpha(n)$, so properties (S5) and (S6) can be proved in the same way as before.

To prove property (S7), since $a_\alpha > c_\alpha^{\nabla h_\alpha}$ for any $\alpha$, it suffices to show that, whenever $\alpha \lhd \beta$ in $2^\omega$, $\lim_{k\to\infty} \frac{a_\alpha(k)}{d_\beta(k)} = 0$. Let $n$ be the minimal number such that $\alpha(n) < \beta(n)$; by the definition of $d_\beta$, $d_\beta(k) \geq (k+1)\cdot a_\alpha(k)$ for any $k \geq n$, which implies that the sequence of the $\frac{a_\alpha(k)}{d_\beta(k)}$ converges to 0. □

## 5. The Grand Finale

Now that we have defined suitable families, we can put everything together to prove the Main Theorem:

**Theorem 5.1.** *Assume* CH *and let $\langle \kappa_\alpha \mid \alpha \in A\rangle$ be a sequence of infinite cardinals such that $|A| \leq \aleph_1$ and $\kappa_\alpha^\omega = \kappa_\alpha$ for every $\alpha \in A$. Given bounding sequences $n_k^-, n_k^+$ and a suitable family $\mathcal{F} = \langle (a_\alpha, d_\alpha, b_\alpha, g_\alpha, f_\alpha, c_\alpha, h_\alpha) \mid \alpha \in A\rangle$, the poset*

$$\mathbb{Q} := \prod_{\alpha\in A} (\mathbb{Q}_{c_\alpha,h_\alpha}^{d_\alpha})^{\kappa_\alpha}$$

*(where the product and all powers have countable support) forces*

$$\mathfrak{v}^\exists_{c_\alpha,h_\alpha} = \mathrm{non}(\mathcal{I}_{f_\alpha}) = \mathfrak{v}^\exists_{b_\alpha,g_\alpha} = \mathfrak{c}^\forall_{a_\alpha,d_\alpha} = \kappa_\alpha$$

*for every $\alpha \in A$.*

(Recall Definition 4.1 for the definition of a suitable family and Definition 3.9 for the definition of the poset $\mathbb{Q}_{c_\alpha,h_\alpha}^{d_\alpha}$.)

This section is dedicated to proving the theorem above. Fix a disjoint family $\langle K_\alpha \mid \alpha \in A\rangle$ such that $|K_\alpha| = \kappa_\alpha$ and let $K := \bigcup_{\alpha\in A} K_\alpha$. Hence, we can express $\mathbb{Q}$ as the countable support product of $\langle \mathbb{Q}_\xi \mid \xi \in K\rangle$, where $\mathbb{Q}_\xi := \mathbb{Q}_{c_\alpha,h_\alpha}^{d_\alpha}$ whenever $\xi \in K_\alpha$. For any $p \in \mathbb{Q}$, $\xi \in \mathrm{supp}(p)$ and $k < \omega$ we write $p(\xi, k) := p(\xi)(k)$.

By the following lemma, without loss of generality we may assume that any $p \in \mathbb{Q}$ is *modest* (unless we explicitly say the opposite), that is, for every $\ell < \omega$ there is at most one $\xi \in \mathrm{supp}(p)$ such that $\ell \in \mathrm{split}(p(\xi))$.

**Lemma 5.2** (modesty)**.** *The set of modest conditions in $\mathbb{Q}$ is dense.*

*Proof.* Fix $p \in \mathbb{Q}$. By induction on $n < \omega$ it is possible to define a function $f = (f_0, f_1)\colon \omega \to \mathrm{supp}(p) \times \omega$ such that

(i) $f_1$ is strictly increasing,
(ii) $f_1(n) \in \mathrm{split}(p(f_0(n)))$ and the norm of $p(f_0(n))$ (with respect to $\mathbb{Q}_{f_0(n)}$) is larger than $n+1$, and
(iii) for any $\xi \in \mathrm{supp}(p)$, the set $f_0^{-1}[\{\xi\}]$ is infinite.

Define $q$ with $\mathrm{supp}(q) = \mathrm{supp}(p)$ such that $q(\xi, k) = p(\xi, k)$ whenever $(\xi, k) \in \mathrm{ran}\, f$, otherwise $q(\xi, k)$ is a singleton contained in $p(\xi, k)$. By the definition of $f$ it is clear that $q \in \mathbb{Q}$ and $q \leq p$. Moreover, by (i), $q$ is modest. □

Thanks to modesty, $\mathbb{Q}$ behaves very much like a single forcing as in section 3. Before we show Theorem 5.1, we revisit the notation and results of section 3 and revise them for the context of $\mathbb{Q}$. The proofs of properness, $\omega^\omega$-bounding and timely reading (as well as the related lemmata) are quite similar to the proofs for the single forcing.

**Notation 5.3.** In the context of $\mathbb{Q}$, for any $p \in \mathbb{Q}$ and $k, n < \omega$:

(1) $\operatorname{split}(p) := \bigcup_{\xi \in \operatorname{supp}(\xi)} \operatorname{split}(p_\xi)$, which is a disjoint union by modesty.
(2) $s_n(p)$ denotes the $n$-th member of $\operatorname{split}(p)$.
(3) For $k \in \operatorname{split}(p)$, let $\xi^p(k)$ denote the unique $\xi \in \operatorname{supp}(p)$ such that $k \in \operatorname{split}(p(\xi))$.
(4) $p{\restriction}^w := \langle p(\xi){\restriction}_w \mid \xi \in \operatorname{supp}(p)\rangle$ when $w$ is a subset of $\omega$ (usually an interval).
(5) $\operatorname{poss}(p, \leq k) := \prod_{\xi \in \operatorname{supp}(p)} \operatorname{poss}(p(\xi), \leq k)$. By modesty, this set is finite with size $\leq \prod_{\ell \leq k} n^+_\ell < n^-_{k+1}$. Define $\operatorname{poss}(p, < k)$ similarly.
(6) For $\eta \in \operatorname{poss}(p, \leq k)$, $p \wedge \eta$ denotes the condition $\langle p(\xi) \wedge \eta(\xi) \mid \xi \in \operatorname{supp}(p)\rangle$.
(7) For $\xi \in K$, $\dot{S}_\xi$ denotes the slalom added by $\mathbb{Q}_\xi$. Note that if $\xi \in K_\alpha$, then $\dot{S}_\xi$ is a $\mathbb{Q}_\xi$-name for a slalom in $\mathcal{S}(c_\alpha, h_\alpha)$.
(8) If $q \in \mathbb{Q}$, $q \leq_n p$ means that $q \leq p$, $\xi^q(s_\ell(q)) \in \operatorname{supp}(p)$ for every $\ell \leq n$, and $(q{\restriction}^{s_n(q)+1}){\restriction}_{\operatorname{supp}(p)} = p{\restriction}^{s_n(q)+1}$.
(9) If $q \in \mathbb{Q}$ and $F \subseteq K$ is finite, $q \leq_{n,F} p$ means that $q \leq p$, $F \subseteq \operatorname{supp}(p)$ and $q(\xi) \leq_n p(\xi)$ for any $\xi \in F$.

**Lemma 5.4.** *If $n < \omega$, $p \in \mathbb{Q}$, and $D \subseteq \mathbb{Q}$ is open dense, then there is some $q \leq_n p$ in $\mathbb{Q}$ such that for any $\eta \in \operatorname{poss}(q, \leq s_n(q))$, $q \wedge \eta \in D$.*

*Proof.* The same argument as in the proof of Lemma 3.4 works here. Concretely, enumerate $\operatorname{poss}(p, \leq s_n(p)) =: \{\eta_k \mid k < m\}$ and, by recursion on $k \leq m$, choose $p_{k+1} \leq \eta'_k{}^\frown p_k{\restriction}^{[s_n(p)+1,\omega)}$ ('pointwise' concatenation) in $D$ such that $s_0(p_{k+1}) > s_n(p)$, $\eta'_k{\restriction}_{\operatorname{supp}(p)} = \eta_k$ and $\eta'_k(\xi)$ is the unique member of $\operatorname{poss}(p_k(\xi), \leq s_n(p))$ for any $\xi \in \operatorname{supp}(p_k) \smallsetminus \operatorname{supp}(p)$. Let $q := r^\frown p_m{\restriction}^{[s_n(p)+1,\omega)}$, where $r$ has domain $\operatorname{supp}(p_m) \times (s_n(p)+1)$, $r(\xi, \ell) := p(\xi, \ell)$ whenever $\xi \in \operatorname{supp}(p)$ and $r(\xi, \ell) := p_m(\xi, \ell)$ otherwise (which is a singleton). □

Define $\mathbb{Q}^*$ as the set of conditions in $p \in \mathbb{Q}$ such that for any $n < \omega$, the norm of $p(\xi^p(s_n(p)), s_n(p))$ is above $n + 1$. Note that the condition $q$ found in Lemma 5.2 is actually in $\mathbb{Q}^*$, so this set is dense in $\mathbb{Q}$. Also, if $p \in \mathbb{Q}^*$, $\alpha \in A$, and $\xi \in K_\alpha \cap \operatorname{supp}(p)$, then $p(\xi) \in \mathbb{Q}^*_\xi := \mathbb{Q}^{*d_\alpha}_{c_\alpha,h_\alpha}$.

**Lemma 5.5** (fusion)**.** *Let $\langle p_n, F_n \mid n < \omega\rangle$ be a sequence such that*

*(i)* $p_n \in \mathbb{Q}^*$ *and* $F_n \subseteq K$ *is non-empty finite,*
*(ii)* $p_{n+1} \leq_{n,F_n} p_n$,
*(iii)* $F_n \subseteq F_{n+1}$ *and* $W := \bigcup_{n<\omega} F_n$ *is equal to* $\bigcup_{n<\omega} \operatorname{supp}(p_n)$.

*Then there is a condition $q \in \mathbb{Q}$ with $\operatorname{supp}(q) = W$ such that $q \leq_{n,F_n} p$ for every $n < \omega$.*

*Proof.* For each $\xi \in W$, let $n_\xi := \min\{n < \omega \mid \xi \in F_n\}$. It is clear that $\langle p_n(\xi) \mid n \geq n_\xi\rangle$ is a sequence in $\mathbb{Q}^*_\xi$ and that $p_{n+1}(\xi) \leq_n p_n(\xi)$. Thus, there is a fusion $q(\xi) \in \mathbb{Q}^*_\xi$ of that sequence as in the proof of Lemma 3.5. Note that $s_k(q(\xi)) = s_k(p_n(\xi))$ for any $n \geq n_\xi$ and $k \leq n$.

Let $q := \langle q(\xi) \mid \xi \in W\rangle$. If $\xi, \zeta \in W$ are different and $k, k' < \omega$, then, for some sufficiently large $n$, we have that $s_k(q(\xi)) = s_k(p_n(\xi)) \neq s_{k'}(p_n(\zeta)) = s_{k'}(q(\zeta))$ (by modesty). Hence $q$ is modest, so $q \in \mathbb{Q}$. □

**Lemma 5.6.** *Let $\chi$ be a sufficiently large regular cardinal and let $N \preceq H_\chi$ be countable such that $\mathbb{Q} \in N$. Let $\langle D_n \mid n < \omega \rangle$ be a sequence of open dense subsets of $\mathbb{Q}$ such that each $D_n \in N$. If $p \in \mathbb{Q} \cap N$, then there is a condition $q \leq p$ in $\mathbb{Q}$ such that for any $n < \omega$ and $\eta \in \operatorname{poss}(q, s_n(q))$, $q \wedge \eta$ is stronger than some member of $D_n \cap N$.*

*Proof.* Enumerate $K \cap N =: \{\xi_k \mid k < \omega\}$ and let $F_n := \{\xi_k \mid k \leq n\}$ for every $n < \omega$. By recursion on $n < \omega$, construct a sequence $\langle p_n \mid n < \omega \rangle$ such that

(i) $p_0 \leq p$,
(ii) $p_n \in \mathbb{Q}^* \cap N$ and $F_n \subseteq \operatorname{supp}(p_n)$,
(iii) $L_n := \langle s_k(p_n(\xi)) \mid k \leq n, \xi \in F_n \rangle$ is an initial segment of $\operatorname{split}(p_n)$,
(iv) for any $k < |L_n|$ and $\eta \in \operatorname{poss}(p_n, {\leq} s_k(p_n))$, $p_n \wedge \eta \in D_k$, and
(v) $p_{n+1} \leq_{n,F_n} p_n$.

Each step of this construction takes place inside $N$, though it is very likely that the final sequence is outside $N$.

Figure 4 illustrates the idea of the following construction. Choose some $p' \leq p$ in $\mathbb{Q}^*$ such that $\xi_0 \in \operatorname{supp}(p')$ and $s_0(p') = s_0(p'(\xi_0))$, and find $p_0 \leq_0 p'$ by application of Lemma 5.4 to $D_0$. Assume that we have constructed $p_n$. Let $m_n := |L_n| - 1 = (n+1)^2 - 1$. By (iv), $p_n$ satisfy the lemma for the splitting levels $k \leq m_n$. The idea in the construction of $p_{n+1}$ is as follows:

(1) We guarantee that $s_k(p_{n+1}) = s_k(p_n)$ for each $k \leq m_n$.
(2) For each $k \leq n$, fix a sufficiently high splitting level of $p_n(\xi_k)$, which will be $s_{m_n+k+1}(p_{n+1})$, and use Lemma 5.4 to guarantee that (iv) holds at the splitting level $m_n + k + 1$.
(3) Add $\xi_{n+1}$ to the support of the condition constructed so far, extend the $\xi_{n+1}$-coordinate so that its first splitting level is sufficiently high (this will be $s_{m_n+n+2}(p_{n+1})$), and use Lemma 5.4 to guarantee that (iv) is satisfied at this level.
(4) At $\xi_{n+1}$, by a recursion of $n+1$ steps choose a sufficiently high splitting level and ensure that (iv) will be satisfied at this level by applying Lemma 5.4. The last splitting level will be $s_{m_n+2(n+1)+1}$.

Note that $m_n + 2(n+1) + 1 = (n+2)^2 - 1$, which is the expected value of $m_{n+1} := |L_{n+1}| - 1$.

We give a formal specification of the construction (with corresponding enumeration). By recursion on $k < 2(n+2)$, construct $p_{n,k}$ as follows:

(1) Let $p_{n,0} := p_n$.
(2) Having defined $p_{n,k}$ for $k \leq n$, choose an $\ell \in \operatorname{split}(p_{n,k}(\xi_k))$ larger than $s_{m_n+k}(p_{n,k})$, choose $p'_{n,k} \leq_{m_n+k} p_{n,k}$ such that $p'_{n,k}(\xi_k, \ell) = p_{n,k}(\xi_k, \ell)$ and $\ell = s_{m_n+k+1}(p'_{n,k})$, and find $p_{n,k+1} \leq_{m_n+k+1} p'_{n,k}$ by application of Lemma 5.4 to $D_{m_n+k+1}$.
(3) Having defined $p_{n,n+1}$, choose $p'_{n,n+1} \leq_{m_n+n+1} p_{n,n+1}$ such that $\xi_{n+1} \in \operatorname{supp}(p'_{n,n+1})$ and $s_0(p'_{n,n+1}(\xi_{n+1})) = s_{m_n+n+2}(p'_{n,n+1})$, and afterwards find $p_{n,n+2} \leq_{m_n+n+2} p'_{n,n+1}$ by application of Lemma 5.4 to $D_{m_n+n+2}$.
(4) Having defined $p_{n,k}$ at $n+2 \leq k < 2(n+1)+1$, choose an $\ell \in \operatorname{split}(p_{n,k}(\xi_{n+1}))$ larger than $s_{m_n+k}(p_{n,k})$, choose $p'_{n,k} \leq_{m_n+k} p_{n,k}$ such that $p'_{n,k}(\xi_{n+1}, \ell) =$

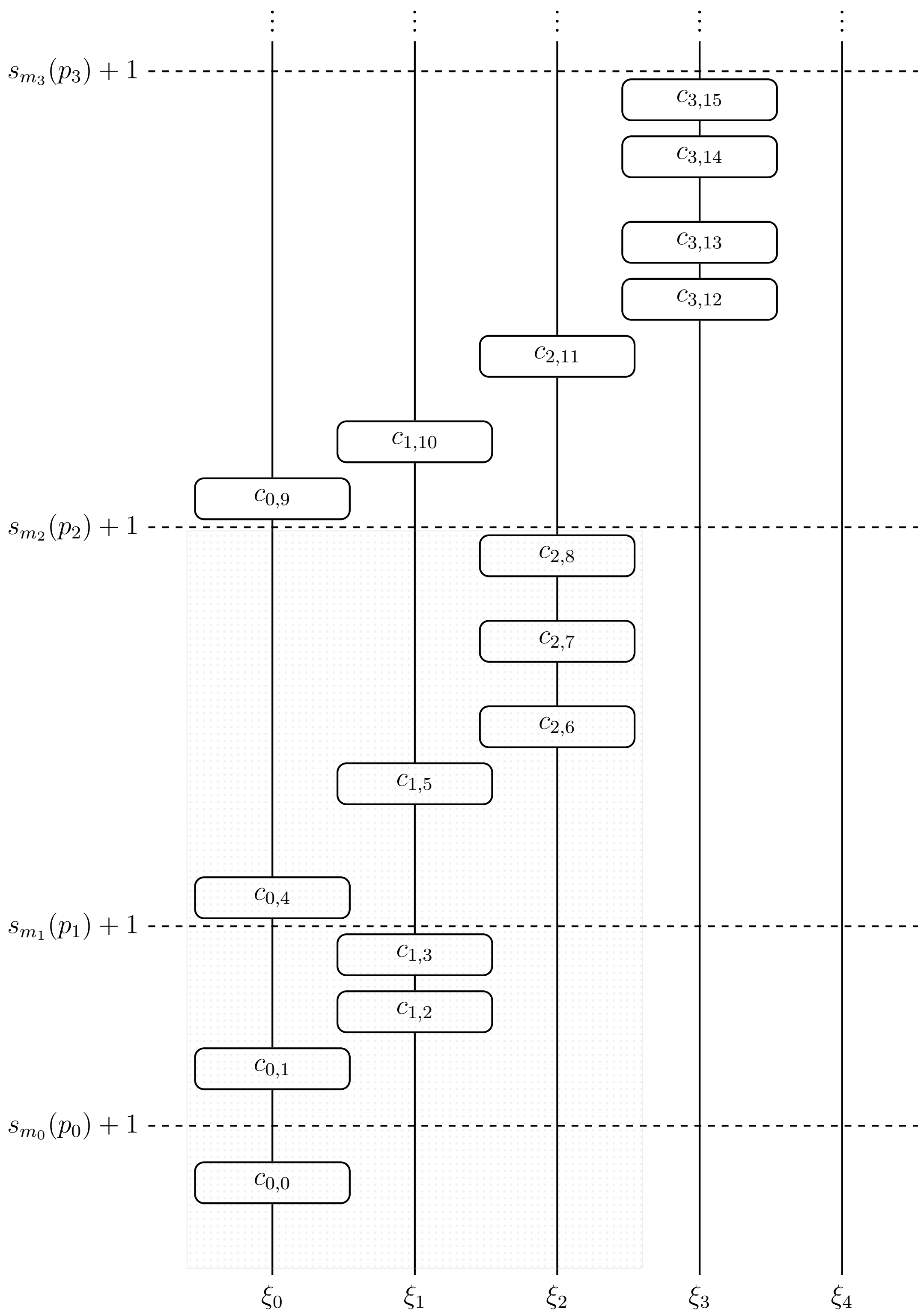

FIGURE 4. The proof of Lemma 5.6 up to the construction of $p_3$. The term "$c_{k,\ell}$" is short for the creature $q(\xi_k, s_\ell(q))$. The shaded area shows the part of $p_2$ that is not modified in the construction of $p_3$, and will also remain constant in any further steps of the construction. Note that $m_0 = 0$, $m_1 = 3$, $m_2 = 8$ and $m_3 = 15$.

$p_{n,k}(\xi_{n+1}, \ell)$ and $\ell = s_{m_n+k+1}(p'_{n,k})$, and find $p_{n,k+1} \leq_{m_n+k+1} p'_{n,k}$ by application of Lemma 5.4 to $D_{m_n+k+1}$.

(5) At the end, $p_{n+1} := p_{n,2(n+1)+1}$ is as desired.

By Lemma 5.5, we can find a fusion $q$ of $\langle p_n \mid n < \omega \rangle$ with $\operatorname{supp}(q) = K \cap N$ such that $q \leq_{n,F_n} p_n$ for any $n < \omega$. By (iii), $q \leq_{|L_n|-1} p_n$ for each $n < \omega$, so (iv) implies that $q$ is as desired. □

**Corollary 5.7.** *$\mathbb{Q}$ is proper and $\omega^\omega$-bounding.*

*Proof.* Properness is immediate by applying Lemma 5.6 to the enumeration of all the open dense subsets in $N$.

To show that $\mathbb{Q}$ is $\omega^\omega$-bounding, let $\dot{x}$ be a $\mathbb{Q}$-name for a member of $\omega^\omega$ and let $p \in \mathbb{Q}$. Let $\chi$ be a sufficiently large regular cardinal and let $N \preceq H_\chi$ be countable such that $\mathbb{Q}, p, \dot{x} \in N$. Let $D_n$ be the set of conditions of $\mathbb{Q}$ that decide $\dot{x}(n)$, which is an open dense set that belongs to $N$. Find $q \leq p$ by application of Lemma 5.6 to $\langle D_n \mid n < \omega \rangle$, which implies that $q \Vdash \dot{x} \in B_n$, where $B_n := \{k < \omega \mid \exists \eta \in \operatorname{poss}(q, s_n(q))\colon q \wedge \eta \Vdash \dot{x}(n) = k\}$ is a finite set. Hence, $q$ forces that $\dot{x}$ is bounded by $f \in V$, where $f(n) := \max(B_n)$. □

**Corollary 5.8** (timely reading of names)**.** *Let $p \in \mathbb{Q}$ and let $\dot{\tau}$ be a $\mathbb{Q}$-name for a function from $\omega$ into the ground model $V$. Then there is a condition $q \leq p$ in $\mathbb{Q}$ such that $q$ reads $\dot{\tau}$ timely.*

*Proof.* A fusion argument as in the proof of Lemma 5.6 works, but here the dense sets need to be defined within the construction of the fusion sequence. Concretely, start with a sufficiently large regular cardinal $\chi$ and a countable model $N \preceq H_\chi$ such that $\mathbb{Q}, p, \dot{\tau} \in N$, and construct the fusion sequence $\langle p_n \mid n < \omega \rangle$ satisfying conditions (i)–(iii) and (v) in the proof of Lemma 5.6 – but instead of (iv) demand

(iv') for any $k < |L_n|$ and $\eta \in \operatorname{poss}(p_n, s_k(p_n))$, $p_n \wedge \eta \in D_k$, where $D_k$ is the (open dense) set of conditions in $\mathbb{Q}^*$ that decide $\dot{\tau}\restriction_{s_k(p_n)}$.

The construction is the same as in Lemma 5.6 except that the dense sets $D_k$ are defined just before each application of Lemma 5.4. □

**Lemma 5.9.** *The poset $\mathbb{Q}$ is $\aleph_2$-cc.*

*Proof.* This follows by a typical $\Delta$-system argument under CH. □

**Lemma 5.10** (few possibilities)**.** *Let $p \in \mathbb{Q}$ and let $\langle g_\xi \mid \xi \in \operatorname{supp}(p) \rangle$ be a sequence of functions from $\omega$ into $\omega$ that go to infinity. Then there is a condition $q \leq p$ in $\mathbb{Q}^*$ such that $\operatorname{supp}(q) = \operatorname{supp}(p)$ and for any $\xi \in \operatorname{supp}(q)$ and $\ell \in \operatorname{split}(q(\xi))$, $|\operatorname{poss}(q, {<}\ell)| < g_\xi(\ell)$.*

*Proof.* Enumerate $\operatorname{supp}(p) =: \{\xi_n \mid n < \omega\}$ and let $F_n := \{\xi_k \mid k \leq n\}$ for every $n < \omega$. By recursion on $n < \omega$, construct a sequence $\langle p_n \mid n < \omega \rangle$ such that

(i) $p_0 \leq p$,

(ii) $p_n \in \mathbb{Q}^*$ and $\operatorname{supp}(p_n) = \operatorname{supp}(p)$,

(iii) $L_n := \langle s_k(p_n(\xi)) \mid k \leq n, \xi \in F_n \rangle$ is an initial segment of $\operatorname{split}(p_n)$,

(iv) for any $\ell \in L_n$, $|\operatorname{poss}(p_n, {<}\ell)| < g_{\xi^{p_n}(\ell)}(\ell)$, and

(v) $p_{n+1} \leq_{n,F_n} p_n$.

Choose $\ell_0 \in \mathrm{split}(p(\xi_0))$ such that both $g_{\xi_0}(\ell_0)$ and the norm of $p(\xi_0, \ell_0)$ are larger than 1, and choose $p_0 \leq p$ in $\mathbb{Q}^*$ with $\mathrm{supp}(p_0) = \mathrm{supp}(p)$ such that $p_0(\xi, k)$ is a singleton for all $k < \ell_0$ and $p_0(\xi_0, \ell_0) = p(\xi_0, \ell_0)$ (so $s_0(p_0) = \ell_0$). Having constructed $p_n$, we can define by recursion an increasing sequence $\langle \ell_{n,k} \mid k < 2n{+}3 \rangle$ of natural numbers such that

(I) $\ell_{n,0} > \max(L_n)$,
(II) $\ell_{n,k} \in \mathrm{split}(p_n(\xi_{k_*}))$ with $k_* := \min\{k, n+1\}$, and
(III) $g_{\xi_{k_*}}(\ell_{n,k})$ is larger than $\mathrm{poss}(p_n, {\leq}\max(L_n)) \cdot \prod_{i<k} |p_n(\xi_{i_*}, \ell_{n,i})|$.

Define $p_{n+1}$ with $\mathrm{supp}(p_{n+1}) = \mathrm{supp}(p_n)$ such that $p_{n+1}(\xi, k) = p_n(\xi, k)$ when either $\xi \in \mathrm{supp}(p_n)$ and $k \in [0, \max(L_n)] \cup [\ell_{n,2n+2}, \omega)$, or $\xi = \xi_{i_*}$ and $k = \ell_{n,i}$ for some $i < 2(n+2)$, and otherwise such that $p_{n+1}(\xi, k)$ is a singleton contained in $p_n(\xi, k)$. It is clear that $p_{n+1}$ is as required.

Define $q$ as in Lemma 5.5. ☐

**Lemma 5.11** (early reading)**.** *Let $\langle X_k \mid k < \omega \rangle$ be a sequence of non-empty sets with $|X_k| \leq n_k^+$, and let $\dot\tau$ be a $\mathbb{Q}$-name for a member of $\prod_{k<\omega} X_k$. Then for any $p \in \mathbb{Q}$, there is a condition $q \leq p$ in $\mathbb{Q}$ which reads $\dot\tau$ early, that is, for any $n < \omega$ and $\eta \in \mathrm{poss}(q, {<}n)$, $q \wedge \eta$ already decides $\dot\tau{\restriction}_n$.*

*Proof.* Without loss of generality, by Corollary 5.8 assume that $p$ reads $\dot\tau$ timely. By recursion on $k < \omega$ we define $q(\xi, k)$ for all $\xi \in \mathrm{supp}(p)$ (so at the end, $\mathrm{supp}(q) = \mathrm{supp}(p)$). When $k \notin \mathrm{split}(p)$, let $q(\xi, k) := p(\xi, k)$. Assume $k \in \mathrm{split}(p)$ and let $p_k := q{\restriction}^k {}^\frown p{\restriction}^{[k,\omega)}$ ('pointwise' concatenation). Exactly as in the proof of Lemma 3.12 (bigness can be used because $|\prod_{i<k} X_i| \leq \prod_{i<k} n_i^+ < n_k^- < d_\alpha(k)$ for any $\alpha \in A$), we can find $q(\xi^p(k), k) \subseteq p(\xi^p(k), k)$ such that the difference between their norms is less that 1 and, for every $\eta \in \mathrm{poss}(p_k, {<}k)$, there is an $r_\eta$ such that $p_k \wedge \eta' \Vdash \dot\tau{\restriction}_k = r_\eta$ for any $\eta' \in \mathrm{poss}(p_k, {\leq}k)$ that extends $\eta$ with $\eta'(\xi^p(k), k) \subseteq q(\xi^p(k), k)$. For $\xi \neq \xi^p(k)$, just define $q(\xi, k) := p(\xi, k)$. ☐

*Proof of Theorem 5.1.* First, we remark that since $\mathbb{Q}$ has $\aleph_2$-cc, it preserves cardinalities and cofinalities. Now note that since $\mathcal{F}$ is suitable, the assumptions in Lemma 2.8 are fulfilled for any $\alpha \in A$: (S1) and (S3) imply (I1); (S4) implies (I2); (S1) and (S5) imply (I3); and (S1) and (S6) imply (I4). Hence $\mathfrak{v}^\exists_{c_\alpha,h_\alpha} \leq \mathrm{non}(\mathcal{I}_{f_\alpha}) \leq \mathfrak{v}^\exists_{b_\alpha,g_\alpha} \leq \mathfrak{c}^\forall_{a_\alpha,d_\alpha}$, so it suffices to show $\mathbb{Q} \Vdash \mathfrak{v}^\exists_{c_\alpha,h_\alpha} \geq \kappa_\alpha$ and $\mathbb{Q} \Vdash \mathfrak{c}^\forall_{a_\alpha,d_\alpha} \leq \kappa_\alpha$ for any $\alpha \in A$.

**The lower bound.** Fix some $\alpha \in A$. We prove $\mathbb{Q} \Vdash \mathfrak{v}^\exists_{c_\alpha,h_\alpha} \geq \kappa_\alpha$ when $\kappa_\alpha > \aleph_1$ (otherwise it is trivial) as follows (analogously to the proof of Lemma 3.3): Let $\dot F$ be a $\mathbb{Q}$-name for a subset of $\prod c_\alpha$ of size $< \kappa_\alpha$. Without loss of generality, we may assume that there is some $\nu < \kappa_\alpha$ such that $\mathbb{Q}$ forces $|\dot F| \leq \nu$.[10] By $\aleph_2$-cc, there is some $B \in V$ with $|B| \leq \max\{\aleph_1, \nu\} < \kappa_\alpha$ such that $\dot F$ is a $\mathbb{Q}{\restriction}_B$-name. Hence there is some $\xi \in K_\alpha \setminus B$.

It suffices to show that for any $\mathbb{Q}{\restriction}_B$-name $\dot x$ for a member of $\prod c_\alpha$, $\mathbb{Q}$ forces that $\dot x(k) \in \dot S_\xi(k)$ for infinitely many $k < \omega$. Only for this argument, we can briefly forget about modesty. Fix $p \in \mathbb{Q}$ and some $n_0 < \omega$ and, without loss of generality, also assume $\xi \in \mathrm{supp}(p)$. Pick a $k \geq n_0$ such that the norm of $p(\xi,k)$ is at least 1; afterwards, find an $r \leq p{\restriction}_B$ which decides $\dot{x}(k) = z$. Define $p' \in \mathbb{Q}$ such that $\mathrm{supp}(p') = \mathrm{supp}(r) \cup \mathrm{supp}(p)$, $p'(\gamma) := r(\gamma)$ for $\gamma \in B$ and $p'(\gamma) := p(\gamma)$ otherwise. Finally, strengthen $p'$ to $q$ with the same support by first setting $q{\restriction}_{K\smallsetminus\{\xi\}} := p'{\restriction}_{K\smallsetminus\{\xi\}}$, and then setting $q(\xi,k) := \{t\}$ for some $t \in p'(\xi,k)$ which contains $z$ (which exists by the definition of the norm, as explained in Lemma 3.3) and $q(\xi,n) := p'(\xi,n)$ when $n \neq k$. Hence $q \leq p'$ and $q \Vdash \dot{x}(k) \in \dot{S}_\xi(k)$.

[10] Concretely, there are $p_0 \in \mathbb{Q}$ and $\nu < \kappa_\alpha$ such that $p_0 \Vdash |\dot F| = \nu$. So we just replace $\dot F$ by some $\mathbb{Q}$-name $\dot F'$ for a subset of $\prod c_\alpha$ of size $\leq \nu$ such that $p_0 \Vdash \dot F' = \dot F$.

This argument shows that given any family of size less than $\kappa_\alpha$, this set cannot be anti-localising in the sense of $\mathfrak{v}^{\exists}_{c_\alpha,h_\alpha}$, and hence $\mathbb{Q} \Vdash \mathfrak{v}^{\exists}_{c_\alpha,h_\alpha} \geq \kappa_\alpha$.

**The upper bound.** This argument is very similar to Lemma 3.13. Fix $\alpha \in A$ and let $C^\alpha := \bigcup\{K_\beta \mid \kappa_\beta \leq \kappa_\alpha,\ \beta \in A\}$. Note that $|C^\alpha| = \kappa_\alpha$ and $\mathbb{Q}{\restriction}_{C^\alpha}$ forces that $\mathfrak{c} = \kappa_\alpha$. To show $\mathbb{Q} \Vdash \mathfrak{c}^{\forall}_{a_\alpha,d_\alpha} \leq \kappa_\alpha$ it suffices to prove that for any $p \in \mathbb{Q}$ and any $\mathbb{Q}$-name $\dot{x}$ for a member of $\prod a_\alpha$, there is some $q \leq p$ in $\mathbb{Q}$ and some $\mathbb{Q}{\restriction}_{C^\alpha}$-name $\dot{\varphi}$ such that $q \Vdash$ "$\dot{\varphi} \in \mathcal{S}(a_\alpha,d_\alpha)$ and $\dot{x} \in^* \dot{\varphi}$". (This means that whenever $G$ is $\mathbb{Q}$-generic over $V$, any member of $\prod a_\alpha$ is localised in $V[G]$ by some slalom in $\mathcal{S}(a_\alpha,d_\alpha) \cap V[G \cap \mathbb{Q}{\restriction}_{C^\alpha}]$, where the latter set has size $\kappa_\alpha$.) Without loss of generality, we may assume that

(i) $p \in \mathbb{Q}^*$ reads $\dot{x}$ early (by Lemma 5.11 because $a_\alpha(k) \leq n_k^+$), and

(ii) for any $\xi \in \mathrm{supp}(p) \smallsetminus C^\alpha$, if $\beta \in A$ and $\xi \in K_\beta$, then for any $k \in \mathrm{split}(p(\xi))$,

$$\min\left\{\frac{c_\beta^{\nabla h_\beta}(k)}{d_\alpha(k)}, \frac{a_\alpha(k)}{d_\beta(k)}\right\} \leq \frac{1}{2\,|\,\mathrm{poss}(p,{<}k)|}$$

(by (S7) and Lemma 5.10 applied to the function $g_\xi$ defined by

$$g_\xi(k) := \left(2 \cdot \min\left\{\frac{c_\beta^{\nabla h_\beta}(k)}{d_\alpha(k)}, \frac{a_\alpha(k)}{d_\beta(k)}\right\}\right)^{-1}$$

when $\xi$ is as above, or the identity on $\omega$ otherwise).

By recursion on $k$, we construct $\dot{\varphi}(k)$ and $q(\xi,k)$ for any $\xi \in \mathrm{supp}(p)$ (so at the end, $\mathrm{supp}(q) = \mathrm{supp}(p)$). Let $p_k := q{\restriction}^{k\frown} p{\restriction}^{[k,\omega)}$. We distinguish two cases.

<u>When $k \notin \mathrm{split}(p)$.</u> In this case, for any $\eta \in \mathrm{poss}(p_k,{<}k)$, $p_k \wedge \eta$ already decides $\dot{x}(k)$ to be some $r_\eta \in a_\alpha(k)$, so we define $\dot{\varphi}(k)$ (in the ground model) as the set of those $r_\eta$. Note that $|\dot{\varphi}(k)| \leq |\,\mathrm{poss}(p_k,<k)| < n_k^- < d_\alpha(k)$. Let $q(\xi,k) := p(\xi,k)$ for any $\xi \in \mathrm{supp}(p)$.

<u>When $k \in \mathrm{split}(p)$.</u> First, assume that $\xi^p(k) \in C^\alpha$. For any $\eta \in \mathrm{poss}(p_k,{<}k)$ and $s \in p(\xi^p(k),k)$, there is a unique $\eta' \in \mathrm{poss}(p_k,{\leq}k)$ that extends $\eta$ and such that $\eta'(\xi^p(k),k) = \{s\}$, so by (i) there is some $r_{\eta,s} \in a_\alpha(k)$ such that $p_k \wedge \eta' \Vdash \dot{x}(k) = r_{\eta,s}$. Define $\dot{r}_\eta$ as the $\mathbb{Q}{\restriction}_{C^\alpha}$-name for a member of $a_\alpha(k)$ such that any condition $p' \in \mathbb{Q}{\restriction}_{C^\alpha}$ with $p'(\xi^p(k),k) = \{s\}$ for some $s \in p(\xi^p(k),k)$ forces that $\dot{r}_\eta = r_{\eta,s}$. Let $\dot{\varphi}(k) := \{\dot{r}_\eta \mid \eta \in \mathrm{poss}(p_k,{<}k)\}$, which is clearly a $\mathbb{Q}{\restriction}_{C^\alpha}$-name for a set of size $< d_\alpha(k)$. Also, for any $\eta \in \mathrm{poss}(p_k,{<}k)$, $p_k \wedge \eta$ forces $\dot{x} \in \dot{\varphi}(k)$. Let $q(\xi,k) := p(\xi,k)$ for all $\xi \in \mathrm{supp}(p)$.

Now assume that $\xi^p(k) \notin C^\alpha$ and let $\beta \in A$ be such that $\xi^p(k) \in K_\beta$. According to (ii), we distinguish two subcases. If

$$|\operatorname{poss}(p_k, <k)| \leq \frac{d_\alpha(k)}{2c_\beta^{\nabla h_\beta}(k)},$$

then $|\operatorname{poss}(p_k, \leq k)| \leq |\operatorname{poss}(p_k, <k)| \cdot c_\beta^{\nabla h_\beta}(k) < d_\alpha(k)$. Let $\dot{\varphi}(k)$ (in the ground model) be the set of objects in $a_\alpha(k)$ that are decided to be $\dot{x}(k)$ by $p_k \wedge \eta$ for some $\eta \in \operatorname{poss}(p_k, \leq k)$. It is clear that this set has size $< d_\alpha(k)$, so we can define $q(\xi, k) := p(\xi, k)$ for all $\xi \in \operatorname{supp}(p)$.

The other subcase is when $|\operatorname{poss}(p_k, <k)| \leq \frac{d_\beta(k)}{2a_\alpha(k)}$. Exactly as in the proof of Lemma 3.13 (with $a = a_\alpha$, $d = d_\beta$ and $e = d_\alpha$), we can find $q(\xi^p(k), k) \subseteq p(\xi^p(k), k)$ such that the difference between their norms is less than 1, as well as a $\dot{\varphi}(k) \in [a_\alpha(k)]^{\leq d_\alpha(k)}$ (in the ground model) such that $p_{k+1} := q{\restriction}^{k+1} {}^\frown p{\restriction}^{[k+1,\omega)}$ forces that $\dot{x} \in \dot{\varphi}(k)$. For any $\xi \in \operatorname{supp}(p) \smallsetminus \{\xi^p(k)\}$, let $q(\xi, k) := p(\xi, k)$.

In the end, both $q$ and $\varphi$ are as required. $\square$

## 6. Open Questions

Concerning the consistency of infinitely many pairwise different cardinals associated with Yorioka ideals, the following summarises the current open questions.

**Question C.** *Is each of the statements below consistent with* ZFC*?*

- *(Q1) There are continuum many pairwise different cardinal invariants of the form* $\operatorname{cov}(\mathcal{I}_f)$.
- *(Q2) There are continuum many pairwise different cardinal invariants of the form* $\operatorname{non}(\mathcal{I}_f)$.
- *(Q3) There are infinitely many pairwise different cardinal invariants of the form* $\operatorname{add}(\mathcal{I}_f)$.
- *(Q4) There are infinitely many pairwise different cardinal invariants of the form* $\operatorname{cof}(\mathcal{I}_f)$.

It is also interesting to consider the consistency of the conjunction of the statements above.

As mentioned in the introduction, Kamo and Osuga [KO14] used the existence of a weakly inaccessible cardinal to force the statement in (Q1), but its consistency is still unknown assuming only ZFC. We believe that extending the ideas in our construction with lim inf techniques as in [KS12] would work to solve this problem (even simultaneously with (Q2)).

Very little is known about the additivity and cofinality of Yorioka ideals. Even the following question is still open (see also [CM19, section 6]).

**Question D.** *Is it consistent with* ZFC *that there are two Yorioka ideals with different additivity numbers (or cofinality numbers)?*

Any idea to solve this question in the positive could be used to prove the consistency of (Q4) using a lim sup creature construction as in this paper. However, as the additivity numbers of Yorioka ideals are below $\mathfrak{b}$, the typical $\omega^\omega$-bounding creature forcing notions will not work to prove the consistency of (Q3).

As for the localisation and anti-localisation cardinals, the remaining open problems are the following.

**Question E.** *Is each of the statements below consistent with* ZFC*?*

*(Q5) There are continuum many pairwise different cardinal invariants of the form* $\mathfrak{v}^{\forall}_{c,h}$.

*(Q6) There are continuum many pairwise different cardinal invariants of the form* $\mathfrak{v}^{\exists}_{c,h}$.

Brendle and the second author [BM14] used a weakly inaccessible cardinal to force (Q5), but its consistency with respect to ZFC alone is still open.

Institute of Discrete Mathematics and Geometry, TU Wien, Wiedner Hauptstrasse 8–10/104, 1040 Wien, Austria

*Email address*: mail@l17r.eu

*URL*: https://l17r.eu

Faculty of Science, Shizuoka University, 836 Ohya, Suruga-ku, Shizuoka City, Shizuoka Prefecture, 422-8529, Japan

*Email address*: diego.mejia@shizuoka.ac.jp